\documentclass[11pt,leqno]{article}
\usepackage{latexsym,amssymb}
\usepackage{amsmath}

\newtheorem{Theorem}{\bf Theorem}[section]
\newtheorem{Lemma}{\bf Lemma}[section]
\newtheorem{Proposition}{\bf Proposition}[section]
\newtheorem{Corollary}{\bf Corollary}[section]
\newtheorem{Remark}{\bf Remark}[section]
\newtheorem{Example}{\bf Example}[section]
\newtheorem{Definition}{\bf Definition}[section]

\newenvironment{theorem}{\begin{Theorem}$\!\!\!$}{\end{Theorem}}
\newenvironment{lemma}{\begin{Lemma}$\!\!\!$}{\end{Lemma}}

\newenvironment{remark}{\begin{Remark}$\!\!\!$}{\end{Remark}}

\newenvironment{definition}{\begin{Definition}$\!\!\!$}{\end{Definition}}

\numberwithin{equation}{section}

\begin{document}

\title{Asymptotic expansions of the solutions of \\
the Cauchy problem for nonlinear parabolic equations\\
\quad}
\author{Kazuhiro Ishige
\thanks{Supported in part by the Grant-in-Aid
for Scientific Research (B)(No. 23340035), Japan Society for the
Promotion of Science.} \\
Mathematical Institute, Tohoku University,\\
Aoba, Sendai 980-8578, Japan\\
(e-mail address : ishige$@$math.tohoku.ac.jp)\\
\quad\\
and\\
\quad\\
Tatsuki Kawakami\\
Department of Mathematical Sciences, Osaka Prefecture University,\\
Sakai 599-8531, Japan\\
(e-mail address : kawakami@ms.osakafu-u.ac.jp)
}
\date{}
\maketitle
%
%
\begin{abstract}
Let $u$ be a solution of the Cauchy problem for the nonlinear parabolic equation  
$$
\partial_t u=\Delta u+F(x,t,u,\nabla u)
\quad\mbox{in}\quad{\bf R}^N\times(0,\infty),
\quad
u(x,0)=\varphi(x)\quad\mbox{in}\quad{\bf R}^N, 
$$
and assume that the solution $u$ behaves like the Gauss kernel as $t\to\infty$. 
In this paper, under suitable assumptions of the reaction term $F$ and the initial function $\varphi$, 
we establish the method of obtaining 
higher order asymptotic expansions of the solution $u$ as $t\to\infty$. 
This paper is a generalization of our previous paper \cite{IK}, 
and our arguments are applicable to the large class of nonlinear parabolic equations. 
\end{abstract}
%
\section{Introduction}
Let $u$ be a unique solution of the Cauchy problem for the nonlinear parabolic equation 
\begin{equation}
\label{eq:1.1}
\left\{
\begin{array}{ll}
\partial_t u=\Delta u+F(x,t,u,\nabla u)\quad & \mbox{in}\quad{\bf R}^N\times(0,\infty),\vspace{3pt}\\
u(x,0)=\varphi(x)\quad & \mbox{in}\quad{\bf R}^N,
\end{array}
\right.
\end{equation}
where $N\ge 1$, $\partial_t=\partial/\partial t$,
$F\in C({\bf R}^N\times(0,\infty)\times {\bf R}\times{\bf R}^N)$, and 
\begin{equation}
\label{eq:1.2}
\varphi\in L^1_K:=
\left\{\phi\in L^1({\bf R}^N):\,
\int_{{\bf R}^N}(1+|x|)^K|\phi(x)|dx<\infty\right\}
\end{equation}
for some constant $K\ge 0$. 
Let $A>1$ and assume that the solution $u$ satisfies 
\begin{equation}
\tag{$C_A$}
 |F(x,t,u(x,t),\nabla u(x,t))|\le C_*(1+t)^{-A}(|u(x,t)|+(1+t)^{1/2}|\nabla u(x,t)|)
\end{equation}
for almost all $(x,t)\in{\bf R}^N\times(0,\infty)$, where $C_*$ is a constant. 
Then it can be proved that 
\begin{eqnarray*}
 & & u\in{\cal S}:=\biggr\{
v\in L^\infty_{loc}(0,\infty:W^{1,\infty}({\bf R}^N))\,:\\
 & & \qquad\qquad\qquad\qquad
\sup_{t>0}\,t^{N/2}\left[\|v(t)\|_{L^\infty({\bf R}^N)}+t^{1/2}\|\nabla v(t)\|_{L^\infty({\bf R}^N)}\right]<\infty
\biggr\},
\end{eqnarray*}
and the solution $u$ behaves like the Gauss kernel as $t\to\infty$, 
that is, 
\begin{equation}
\label{eq:1.3}
\left\{
\begin{array}{l}
\mbox{$\displaystyle{\int_{{\bf R}^N}u(x,t)dx}$ converges to a constant $M$ as $t\to\infty$ and}\vspace{7pt}\\
\mbox{$\displaystyle{\lim_{t\to\infty}}
\|u(t)-MG(1+t)\|_{L^q({\bf R}^N)}/\|G(1+t)\|_{L^q({\bf R}^N)}=0$ for any $q\in[1,\infty]$},
\end{array}
\right.
\end{equation}
where 
$$
G(x,t)=(4\pi t)^{-\frac{N}{2}}\exp\left(-\frac{|x|^2}{4t}\right)
$$ 
(see Theorem~\ref{Theorem:3.1}). 
We introduce the condition $(F_A)$ on the reaction term $F$:
\begin{equation}
\tag{$F_A$}
\left\{
\begin{array}{ll}
{\rm (i)} &  \mbox{$F(x,t,0,0)=0$ for all $(x,t)\in{\bf R}^N\times(0,\infty)$;}\vspace{5pt}\\
{\rm (ii)} & \mbox{For any $v_1$ and $v_2\in {\cal S}$, there exists a constant $C$ such that}\vspace{5pt}\\
 & \,\,\,
 |F(x,t,v_1(x,t),\nabla v_1(x,t))-F(x,t,v_2(x,t),\nabla v_2(x,t))|\vspace{5pt}\\
 & \quad
 \le C(1+t)^{-A}(|v_1(x,t)-v_2(x,t)|+(1+t)^{1/2}|\nabla v_1(x,t)-\nabla v_2(x,t)|)\vspace{5pt}\\
 & \mbox{for almost all $(x,t)\in{\bf R}^N\times(0,\infty)$}.
\end{array}
\right.
\end{equation}
Condition $(F_A)$ ensures that, 
if $v\in{\cal S}$, then $v$ satisfies condition $(C_A)$. 
In this paper, under these conditions $(C_A)$ and $(F_A)$, 
we study the large time behavior of the solution $u$ of \eqref{eq:1.1}, 
and establish the method of obtaining higher order asymptotic expansions of the solution $u$ as $t\to\infty$. 
\vspace{3pt}

Consider the Cauchy problem for the semilinear heat equation 
\begin{equation}
\label{eq:1.4}
\partial_t u=\Delta u+\lambda|u|^{p-1}u\quad\mbox{in}\quad{\bf R}^N\times(0,\infty),
\quad
u(x,0)=\varphi(x)\quad\mbox{in}\quad{\bf R}^N,
\end{equation}
where $N\ge 1$, $\lambda\in{\bf R}$, $p>1+2/N$, and $\varphi\in L^1({\bf R}^N)\cap L^\infty({\bf R}^N)$. 
Under suitable assumptions, 
Cauchy problem \eqref{eq:1.4} has a unique global in time solution, and 
the large time behavior of the solution has been studied in many papers by various methods 
(see for example \cite{DolbeaultKarch},
\cite{MR0913672}, \cite{Fujita}--\cite{IK}, \cite{MR0837255},
\cite{MR0921547}--\cite{kawa97}, \cite{Qit}--\cite{Ta}, \cite{z}, and references therein). 
In particular, it is known that, 
if 
$$
\varphi\in L^1({\bf R}^N)\cap L^\infty({\bf R}^N)\quad\mbox{and}\quad
\mbox{$\|\varphi\|_{L^{N(p-1)/2}({\bf R}^N)}$ is sufficiently small}, 
$$
then there exists a unique global in time solution of \eqref{eq:1.4}, satisfying \eqref{eq:1.3}. 
In \cite{IIK} 
the authors of this paper and Ishiwata studied 
the large time behavior of the solution of \eqref{eq:1.4},  
and investigated the decay rate of the difference 
between the solution $u$ satisfying \eqref{eq:1.3} and the Gauss kernel 
(see also \cite{IshigeKawakami}, \cite{MR1373470}, \cite{kawa97}, \cite{Ta}, and \cite[Proposition 20.13]{MR2346798}). 
Subsequently, in \cite{IK},  improving the arguments in \cite{IIK}, 
the authors of this paper studied the Cauchy problem for the nonlinear parabolic equations of type 
$$
\partial_t u=\Delta u+F(x,t,u)
\quad\mbox{in}\quad{\bf R}^N\times(0,\infty),
$$
and gave higher order asymptotic expansions of the solution satisfying \eqref{eq:1.3}.  
Their results are applicable to the solution of \eqref{eq:1.4}, satisfying \eqref{eq:1.3}. 
We remark that, if the solution $u$ of \eqref{eq:1.4} satisfies \eqref{eq:1.3}, 
then there holds 
$$
\left|\lambda|u(x,t)|^{p-1}u(x,t)\right|\le C(1+t)^{-\frac{N}{2}(p-1)}|u(x,t)|,\qquad (x,t)\in{\bf R}^N\times(0,\infty)
$$
for some constant $C$, 
and conditions $(C_A)$ and $(F_A)$ are satisfied with $A=N(p-1)/2>1$. 

On the other hand, for the Cauchy problem for the nonlinear parabolic equations of type 
\begin{equation}
\label{eq:1.5}
\partial_t u=\Delta u+\nabla\cdot\mbox{\boldmath $F$}(x,t,u)
\quad\mbox{in}\quad{\bf R}^N\times(0,\infty),
\end{equation}
under suitable assumptions on $\mbox{\boldmath $F$}$ and the initial function, 
there exists a global in time solution satisfying \eqref{eq:1.3}, 
and the asymptotics of the solution has been studied in detail by many mathematicians 
(see for example \cite{Carpio}, \cite{Carpio2}, \cite{DC}, \cite{DZ}, \cite{EVZ}, \cite{EZ}, \cite{FM}, 
\cite{MR2483011}, \cite{NSY}, \cite{NY}, \cite{Yamada}, \cite{Yamada2}, \cite{Zua}, and references therein).   
The solution $u$ of the Cauchy problem for \eqref{eq:1.5} satisfies 
\begin{equation}
\label{eq:1.6}
\int_{{\bf R}^N}u(x,t)dx=\int_{{\bf R}^N}u(x,0)dx
\end{equation}
under suitable integrability conditions on the solution $u$, and property \eqref{eq:1.6} has been used effectively 
in the study of the asymptotic expansions of the solution of \eqref{eq:1.5} in the papers. 
However the solution of \eqref{eq:1.1} does not necessarily have property \eqref{eq:1.6}, 
and it seems difficult to apply their arguments to Cauchy problem \eqref{eq:1.1}
for general nonlinear parabolic equations directly. 
\vspace{3pt}

This paper is a generalization of our previous paper \cite{IK}, 
and the main results of this paper are given in Section~4. 
In this paper, 
by using the operator $P_{[K]}(t)$ introduced by \cite{IIK} (see Section~2.1) 
we establish the method of obtaining higher order asymptotic expansions 
of the solution of Cauchy problem \eqref{eq:1.1} under conditions $(C_A)$ and $(F_A)$. 
Furthermore we give decay estimates of the difference between the solution and its asymptotic expansions. 
Our results can give not only higher order asymptotic expansions of the solutions 
of general nonlinear parabolic equations systematically but also sharp asymptotic expansions  
of the solutions for some typical examples of nonlinear parabolic equations. 
In Section~6 we apply our results to some selected examples of nonlinear parabolic equations 
including the convection-diffusion equation and the Keller-Segel system of parabolic-parabolic type, 
and explain the advantage of our results. 
\vspace{3pt}

The rest of this paper is organized as follows. 
In Section~2 we give some notation and introduce the operator $P_{[K]}(t)$. 
Furthermore we recall some properties of the solution of the heat equation and the operator $P_{[K]}(t)$, 
and give a preliminary lemma on the volume potential (see also Section~7). 
%
%
In Section 3 we give a theorem,
%
%
which implies that the solution of \eqref{eq:1.1} belongs to ${\cal S}$ and satisfies \eqref{eq:1.3} 
and which ensures the well-definedness of $P_{[K]}(t)u(t)$ and $P_{[K]}(t)F(\cdot,t,u(t),\nabla u(t))$. 
In Section~4 we state the main results of this paper, 
and give higher order asymptotic expansions of the solution $u$ of \eqref{eq:1.1} 
under conditions $(C_A)$ and $(F_A)$ with $A>1$.
Section 5 is devoted to the proof of theorems given in Section~4. 
In Section~6 we apply our main results to some selected examples of nonlinear parabolic equations. 
Section~7 is an appendix, and there we prove the H\"older continuity of the gradient of the volume potential. 
\section{Notation and preliminary results}
In this section we give some notation and the definition of the solution of \eqref{eq:1.1}. 
Furthermore we introduce an operator $P_{[K]}(t)$, and 
recall some preliminary lemmas on the solution of the heat equation and the operator $P_{[K]}(t)$. 
\subsection{Notation and operator $P_{[K]}(t)$}
We introduce some notation. 
Let ${\bf N}_0={\bf N}\,\cup\,\{0\}$. 
For any $k\in{\bf R}$, let $[k]$ be an integer such that $k-1<[k]\le k$. 
For any multi-index $\alpha=(\alpha_1,\cdots,\alpha_N)\in{\bf N}_0^N$, 
we put 
$$
\begin{array}{l}
|\alpha|:=\displaystyle{\sum_{i=1}^N}|\alpha_i|,\quad
\alpha!:=\prod_{i=1}^N\alpha_i!,\quad
x^\alpha:=\prod_{i=1}^N x_i^{\alpha_i},\quad
\partial_x^\alpha:=
\frac{\partial^{|\alpha|}}{\partial x_1^{\alpha_1}\cdots\partial x_N^{\alpha_N}},\vspace{5pt}\\
J(\alpha):=\{\rho=(\rho_1,\cdots,\rho_N)\in{\bf N}_0^N\setminus\{\alpha\}:\rho_i\le\alpha_i
\,\,\mbox{for all}\,\,i=1,\cdots,N\},\vspace{5pt}\\
\displaystyle{g_\alpha(x,t):=\frac{(-1)^{|\alpha|}}{\alpha!}(\partial_x^\alpha G)(x,1+t)}.
\end{array}
$$ 
In particular, we write $g(x,t)=g_0(x,t)$ for simplicity. 
We denote by $e^{t\Delta}\varphi$ 
the unique bounded solution of the Cauchy problem for the heat equation with the initial function $\varphi\in L^1({\bf R}^N)$, 
that is, 
\begin{equation}
\label{eq:2.1}
(e^{t\Delta}\varphi)(x):=\int_{{\bf R}^N}G(x-\xi,t)\varphi(\xi)d\xi.
\end{equation}
For any two nonnegative functions 
$f_1$ and $f_2$ defined in a subset $D$ of $[0,\infty)$, 
we say 
$f_1(t)\preceq f_2(t)$ for all $t\in D$ 
if there exists a positive constant $C$ such that 
$f_1(t)\le Cf_2(t)$ for all $t\in D$. 
In addition, we say  
$f_1(t)\asymp f_2(t)$ for all $t\in D$ 
if $f_1(t)\preceq  f_2(t)$ and $f_2(t)\preceq f_1(t)$ for all $t\in D$. 
In what follows, we write 
$$
\|\cdot\|_q=\|\cdot\|_{L^q({\bf R}^N)},
\qquad
|||\cdot|||_m=\|\cdot\|_{L^1({\bf R}^N,(1+|x|)^mdx)}
$$
for simplicity, where $q\in[1,\infty]$ and $m\ge 0$. 
\vspace{5pt}

We give the definition of the solution of Cauchy problem \eqref{eq:1.1}. 
\begin{definition}
\label{Definition:2.1} 
Let $\varphi\in L^1({\bf R}^N)$ and assume $F\in C({\bf R}^N\times(0,\infty)\times{\bf R}\times{\bf R}^N)$. 
Then the function $u\in L^\infty_{loc}(0,\infty:W^{1,1}({\bf R}^N))$ is said to be a solution of \eqref{eq:1.1} if 
$$
u(x,t)=\int_{{\bf R}^N}
G(x-\xi,t)\varphi(\xi)d\xi+\int_0^t\int_{{\bf R}^N}G(x-\xi,t-s)F(\xi,s,u(\xi,s),\nabla u(\xi,s))d\xi ds
$$
holds for almost all $(x,t)\in{\bf R}^N\times(0,\infty)$.
\end{definition}
\vspace{5pt}

Let $k\in{\bf N}_0$, $i\in\{0,\dots,k\}$, and $t>0$.  
Next we follow \cite{IIK} and \cite{IK}, 
and introduce a linear operator $P_i(t)$ on $L^1_k$ by
\begin{equation}
\label{eq:2.2}
[P_i(t)f](x):=f(x)-\sum_{|\alpha|\le i}M_\alpha(f,t) g_\alpha(x,t),
\end{equation}
where $f\in L^1_k$ and $M_\alpha(f,t)$ is the constant defined inductively (in $\alpha$) by
\begin{equation}
\label{eq:2.3}
\begin{array}{l}
\displaystyle{M_0(f,t):=\int_{{\bf R}^N}f(x)dx},\qquad
\displaystyle{M_\alpha(f,t):=\int_{{\bf R}^N}x^\alpha f(x)dx}\quad\mbox{if}\quad|\alpha|=1,\vspace{7pt}\\
\displaystyle{M_\alpha(f,t):=\int_{{\bf R}^N} 
x^\alpha f(x)dx-\sum_{\rho\in J(\alpha)}M_\rho(f,t) 
\int_{{\bf R}^N}x^\alpha g_\rho(x,t)dx}\quad\mbox{if}\quad|\alpha|\ge 2.
\end{array}
\end{equation}
Then the operator $P_i(t)$ has the following property,
\begin{equation}
\label{eq:2.4}
\int_{{\bf R}^N}x^\alpha[P_i(t)f](x) dx=0,
\qquad|\alpha|\le i,
\end{equation}
which is a crucial property in our analysis. 
Here, under the assumption $\varphi\in L^1_K$ with $K\ge 0$, 
we apply the operator $P_{[K]}(t)$ to $e^{t\Delta}\varphi$, and obtain 
\begin{eqnarray*}
P_{[K]}(t)e^{t\Delta}\varphi &\!\!\!=\!\!\!& 
e^{t\Delta}\varphi-\sum_{|\alpha|\le[K]}M_\alpha(e^{t\Delta}\varphi,t) g_\alpha(x,t)\\
 &\!\!\!=\!\!\!& e^{t\Delta}\varphi-\sum_{|\alpha|\le[K]}M_\alpha(\varphi,0) g_\alpha(x,t)
 =e^{t\Delta}[P_{[K]}(0)\varphi]
\end{eqnarray*}
for all $t>0$. (See also Lemma~\ref{Lemma:2.3}~(ii).)
Then, due to property \eqref{eq:2.4}, 
we have 
\begin{equation}
\label{eq:2.5}
t^{\frac{N}{2}(1-\frac{1}{q})}
\biggr\|e^{t\Delta}\varphi-\sum_{|\alpha|\le [K]}M_\alpha(\varphi,0)g_\alpha(t)\biggr\|_q=
\left\{
\begin{array}{ll}
o(t^{-\frac{K}{2}}) & \mbox{if}\quad K=[K],\vspace{3pt}\\
O(t^{-\frac{K}{2}})& \mbox{if}\quad K>[K],
\end{array}
\right.
\end{equation}
as $t\to\infty$. This is easily obtained by Lemma~\ref{Lemma:2.1} and 
property $(G1)$ given in Section~2.2. See also \cite[Proposition 2.1]{IK}. 
\subsection{Preliminaries}
In this section we recall some preliminary results 
on the behavior of solutions for the heat equation and the operator $P_{[K]}(t)$. 
Furthermore we give preliminary lemmas on the volume potential and an integral inequality. 

Let $\alpha\in{\bf N}_0^N$ and $g_\alpha$ be the function given in Section~2.1.
Then, for any $j=0,1,2,\dots$,
there exists a constant $C_1$ such that 
\begin{equation} 
\label{eq:2.6}
|\partial_t^j\partial_x^\alpha G(x,t)|\le C_1t^{-\frac{N+|\alpha|+2j}{2}}
\left[1+\left(\frac{|x|}{t^{1/2}}\right)^{|\alpha|+2j}\right]
\exp\left(-\frac{|x|^2}{4t}\right)
\end{equation} 
for all $(x,t)\in{\bf R}^N\times(0,\infty)$. 
This inequality yields the inequalities
\begin{equation} 
\label{eq:2.7}
\|g_\alpha(t)\|_q\preceq (1+t)^{-\frac{N}{2}(1-\frac{1}{q})-\frac{|\alpha|}{2}},
\quad
\int_{{\bf R}^N}|x|^l|g_\alpha(x,t)|dx\preceq (1+t)^{\frac{l-|\alpha|}{2}},
\quad t>0,
\end{equation} 
for any $q\in[1,\infty]$ and $l\ge 0$. 
Furthermore, by \eqref{eq:2.1} and \eqref{eq:2.6}
we have: 
\begin{itemize}
  \item[$(G1)$] For any multi-index $\alpha$ and $1\le p\le q\le\infty$, 
there exists a constant $c_{|\alpha|}$, independent of $p$ and $q$, such that
$$
\|\partial_x^\alpha e^{t\Delta}\varphi\|_q
\le c_{|\alpha|} t^{-\frac{N}{2}(\frac{1}{p}-\frac{1}{q})-\frac{|\alpha|}{2}}\|\varphi\|_p,\quad t>0.
$$
%
%
In particular, there holds 
$\|e^{t\Delta}\varphi\|_q\le\|\varphi\|_q$ for all $t>0$;
\item[$(G2)$] For any $l\ge 0$ and $\delta>0$, 
there exists a constant $C_2$ such that
$$
\int_{{\bf R}^N}|x|^l|(e^{t\Delta}\varphi)(x)|dx
\le(1+\delta)\int_{{\bf R}^N}|x|^l|\varphi(x)|dx
+C_2t^{\frac{l}{2}}\int_{{\bf R}^N}|\varphi(x)|dx,\quad t>0
$$
(see also Lemma 2.1 in \cite{IIK}).
This inequality implies that
$$
|||e^{t\Delta}\varphi|||_l\le (1+\delta)|||\varphi|||_l+C_3(1+t^{\frac{l}{2}})\|\varphi\|_1,\quad t>0,
$$
for some constant $C_3$;
\item[$(G3)$] For any $l\ge 0$, there exists a constant $C_4$ such that 
$$
\int_{{\bf R}^N}|x|^l|\nabla(e^{t\Delta}\varphi)(x)|dx
\le C_4t^{-\frac{1}{2}}\int_{{\bf R}^N}|x|^l|\varphi(x)|dx
+C_4t^{\frac{l-1}{2}}\int_{{\bf R}^N}|\varphi(x)|dx,\quad t>0. 
$$
This inequality implies that
$$
|||\nabla(e^{t\Delta}\varphi)|||_l
\le C_5t^{-\frac{1}{2}}|||\varphi|||_l+C_5t^{-\frac{1}{2}}(1+t^{\frac{l}{2}})\|\varphi\|_1,\quad t>0,
$$
for some constant $C_5$.
\end{itemize}
Moreover we give one lemma on $e^{t\Delta}\varphi$. 
See \cite[Lemmas 2.2 and 2.5]{IIK}. 
\begin{Lemma}
\label{Lemma:2.1}
Let $\varphi\in L^1_k$ with $k\ge 0$ and assume 
$$
\int_{{\bf R}^N}x^\alpha\varphi(x)dx=0,\qquad |\alpha|\le m, 
$$
for some integer $m\in\{0,\dots,[k]\}$. 
Then there holds the following:
\newline
{\rm (i)} 
If $0\le m\le [k]-1$, 
for any $l\in[0,k-m-1]$, 
there exists a constant $C_1$ such that 
$$
\begin{array}{l}
 \displaystyle
 \int_{{\bf R}^N}|x|^l\left|(e^{t\Delta}\varphi)(x)\right|dx\vspace{5pt}\\
 \displaystyle
 \le C_1t^{-\frac{m+1}{2}}
 \left[\int_{{\bf R}^N}|x|^{m+l+1}|\varphi(x)|dx
 +t^{\frac{l}{2}}\int_{{\bf R}^N}|x|^{m+1}|\varphi(x)|dx\right],
 \quad t>0;
\end{array}
$$
{\rm (ii)} If $m=[k]$, 
for any $l\in[0,k-[k]]$, 
there exists a constant $C_2$ such that 
$$
\int_{{\bf R}^N}|x|^l\left|(e^{t\Delta}\varphi)(x)\right|dx
\le C_2t^{-\frac{k-l}{2}}\int_{{\bf R}^N}|x|^k|\varphi(x)|dx
$$
for all $t>0$. 
In particular, if $k=[k]$, then
$\displaystyle{\lim_{t\to\infty}t^{\frac{k}{2}}}
\|e^{t\Delta}\varphi\|_1=0$. 
\end{Lemma}
Next we recall the following two lemmas on the operator $P_k(t)$. 
See \cite[Lemma~2.3]{IIK} and \cite[Lemma~2.3]{IK}. 
\begin{lemma}
\label{Lemma:2.2}
Let $K\ge 0$ and $f$ be a measurable function in ${\bf R}^N\times(0,\infty)$ such that 
$f(t)\in L^1_K$ for all $t>0$.
Then there holds the following: 
\newline
{\rm (i)} 
Assume that there exist constants $\beta\ge0$ and $\gamma\ge 0$ such that 
$$
\sup_{t>0}\,(1+t)^{-\frac{l}{2}+\gamma}t^{\beta}|||f(t)|||_l<\infty
$$
for all $l\in[0,K]$. 
Then, for any multi-index $\alpha$ with $|\alpha|\le [K]$, 
there exists a constant $C_1$ such that 
\begin{equation*}
|M_\alpha(f(t),t)|\le C_1(1+t)^{\frac{|\alpha|}{2}-\gamma}t^{-\beta},\qquad t>0.
\end{equation*}
Furthermore 
$$
\sup_{t>0}\,
\left[t^{\frac{N}{2}(1-\frac{1}{q})+\gamma+\beta}\|P_{[K]}(t)f(t)-f(t)\|_q
+(1+t)^{-\frac{l}{2}+\gamma}t^{\beta}|||P_{[K]}(t)f(t)|||_l\right]<\infty
$$
for any $l\in[0,K]$ and $q\in[1,\infty]$; 
\vspace{3pt}
\newline
{\rm (ii)} If there exist constants $\beta'\ge0$ and $\gamma'\ge 0$ such that 
$$
\sup_{t>0}\left[t^{\frac{N}{2}(1-\frac{1}{q})+\gamma'+\beta'}\|f(t)\|_q
+(1+t)^{-\frac{l}{2}+\gamma'}t^{\beta'}|||f(t)|||_l\right]<\infty
$$
for all $l\in[0,K]$ and $q\in[1,\infty]$, 
then 
$$
t^{\frac{N}{2}(1-\frac{1}{q})+\frac{j}{2}}\left\|\nabla^j\int_0^t e^{(t-s)\Delta}P_{[K]}(s)f(s)ds\right\|_q
\preceq t^{-\frac{K}{2}}\int_0^t(1+s)^{\frac{K}{2}-\gamma'}s^{-\beta'}ds,\qquad t>0,
$$
for any $q\in[1,\infty]$ and $j=0,1$. 
\end{lemma}
\begin{Lemma}
\label{Lemma:2.3}
Let $k\ge0$ and $f=f(x,t)\in C({\bf R}^N\times(0,\infty))\cap L^\infty({\bf R}^N\times(0,\infty))$ such that 
$\sup_{0<\tau<t}|||f(\tau)|||_k<\infty$ for all $t>0$.
Let $u$ be a solution of the Cauchy problem 
$$
\partial_t u=\Delta u+f\quad\mbox{in}\quad{\bf R}^N\times(0,\infty),
\qquad
u(x,0)=\varphi(x)\quad\mbox{in}\quad{\bf R}^N,
$$
where 
$\varphi\in L^1_k$. 
Then there holds the following: 
\newline
$(\rm{i})$ For any $i\in\{0,\cdots,[k]\}$, 
the function $v=[P_i(t)u(t)](x)$ satisfies 
$$
\partial_t v=\Delta v+P_i(t)f(t)\qquad
\mbox{in}\quad{\bf R}^N\times(0,\infty);
$$
$(\rm{ii})$
For any multi-index $\alpha$ with $|\alpha|\le [k]$,
$$
M_\alpha(u(t),t)-M_\alpha(u(s),s)=\int_s^t M_\alpha(f(\tau),\tau)d\tau
$$
for all $t>s\ge 0$. 
In particular, if $f\equiv 0$, 
$$
M_\alpha(u(t),t)
=M_\alpha(\varphi,0),\qquad|\alpha|\le [k],\quad t>0.
$$
\end{Lemma}

Next we give one lemma on the volume potential.  
Let $T>0$ and $H\in L^\infty(0,T:L^\infty({\bf R}^N))$. 
Let $w$ be the the volume potential of $H$ defined by 
\begin{equation}
\label{eq:2.8}
w(x,t):=\int_0^t\int_{{\bf R}^N} G(x-\xi,t-\tau)H(\xi,\tau)d\xi d\tau,\quad t\in(0,T). 
\end{equation}
Then we have: 
\begin{lemma}
\label{Lemma:2.4} 
Let $T>0$ and $H\in L^\infty(0,T:L^\infty({\bf R}^N))$. 
Then $w$ and $\nabla_xw$ are continuous functions in ${\bf R}^N\times(0,T)$
%
%
and 
\begin{equation}
\label{eq:2.9}
(\nabla_x w)(x,t)=\int_0^t\int_{{\bf R}^N}(\nabla_x G)(x-\xi,t-\tau)H(\xi,\tau)d\xi d\tau
\end{equation}
holds for all $(x,t)\in{\bf R}^N\times(0,T)$. 
Furthermore there exists a constant $C_1$
such that  
\begin{equation}
\label{eq:2.10}
\sup_{0<t<T}\|w(t)\|_\infty+
\sup_{0<t<T}\|(\nabla_x w)(t)\|_\infty\le C_1\|H\|_{L^\infty(0,T:L^\infty({\bf R}^N))}.  
\end{equation}
In addition, 
for any $\nu\in(0,1)$ and $|\alpha|\le 1$, 
there exists a constant $C_2$ such that 
\begin{equation}
\label{eq:2.11}
\frac{|\partial_x^\alpha w(x,t)-\partial_x^\alpha w(y,s)|}{|x-y|^\nu+|t-s|^{\nu/2}}\le C_2\|H\|_{L^\infty(0,T:L^\infty({\bf R}^N))}
\end{equation}
for all $(x,t)$, $(y,s)\in{\bf R}^N\times(0,T)$ with $(x,t)\not=(y,s)$. 
\end{lemma}
Lemma~\ref{Lemma:2.4} is proved by the same argument as in  \cite[Chapter~1]{F}. 
We give the proof in Section~7 for completeness of this paper. 

At the end of this section we recall one lemma on an integral inequality. 
See \cite[Lemma~2.4]{IK}. 
\begin{lemma}
\label{Lemma:2.5}
Let $\zeta$ be a nonnegative function in $(0,\infty)$ such that $\sup_{0<t<1}\zeta(t)<\infty$. 
Let $A>1$ and $\sigma>0$. 
If, for any $\delta>0$, there holds 
$$
\zeta(2t)\le (1+\delta)\zeta(t) +C_1\int_t^{2t}s^{-A}\zeta(s)ds+C_1t^\sigma,\quad t\ge 1/2,
$$
for some constant $C_1$, 
then there exists a constant $C_2$ such that
$\zeta(t)\le C_2t^\sigma$ for all $t\ge 1$. 
\end{lemma}
%
\section{Large time behavior of solutions}
Consider the Cauchy problem 
\begin{equation}
\label{eq:3.1}
\left\{
\begin{array}{ll}
\partial_t u=\Delta u+f(x,t,u,\nabla u)\quad & \mbox{in}\quad{\bf R}^N\times(0,\infty),\vspace{3pt}\\
u(x,0)=\varphi(x)\quad & \mbox{in}\quad{\bf R}^N,
\end{array}
\right.
\end{equation}
where $f\in C({\bf R}^N\times(0,\infty)\times{\bf R}\times{\bf R}^N)$ and $\varphi\in L_K^1$ for some $K\ge 0$. 
In this section 
we assume that there exist constants $C>0$ and $A>1$ such that 
\begin{equation}
\label{eq:3.2}
|f(x,t,p,q)|\le C(1+t)^{-A}(|p|+(1+t)^{1/2}|q|)
\end{equation}
for all $(x,t,p,q)\in{\bf R}^N\times(0,\infty)\times{\bf R}\times{\bf R}^N$, 
and prove the following theorem, which ensures the well-definedness of $P_{[K]}(t)u(t)$ 
and $P_{[K]}(t)F(\cdot,t,u(t),\nabla u(t))$
for the solution $u$ of \eqref{eq:1.1} in Section~4. 
\begin{theorem}
\label{Theorem:3.1}
Assume $\varphi\in L_K^1$ for some $K\ge 0$ and condition \eqref{eq:3.2}. 
Then there exists a solution $u$ of \eqref{eq:3.1} with the following properties: \vspace{3pt}
\newline
{\rm (i)} $u$, $\nabla u\in C({\bf R}^N\times(0,\infty))$;\vspace{3pt}
\newline
{\rm (ii)} For any $q\in[1,\infty]$ and $l\in[0,K]$, there hold 
\begin{eqnarray}
\label{eq:3.3}
 & & \sup_{0<t<\infty}t^{\frac{N}{2}(1-\frac{1}{q})}\left[\|u(t)\|_q+t^{\frac{1}{2}}\|(\nabla_x u)(t)\|_q\right]<\infty,\\
\label{eq:3.4}
 & & \sup_{0<t<\infty}(1+t)^{-\frac{l}{2}}
 \left[|||u(t)|||_l+t^{\frac{1}{2}}|||(\nabla_x u)(t)|||_l\right]<\infty;
\end{eqnarray}
{\rm (iii)} There exists a limit 
$$
M:=\lim_{t\to\infty}\int_{{\bf R}^N}u(x,t)dx
=\int_{{\bf R}^N}\varphi(x)dx+\int_0^\infty\int_{{\bf R}^N}f(x,t,u,\nabla u)dxdt
$$
such that
\begin{equation}
\label{eq:3.5}
\lim_{t\to\infty}t^{\frac{N}{2}(1-\frac{1}{q})+\frac{j}{2}}
\left\|\nabla^j\left[u(t)-Mg(t)\right]\right\|_q=0\quad
\mbox{for any $q\in[1,\infty]$ and $j=0,1$.} 
\end{equation}
\end{theorem}
In order to prove Theorem~\ref{Theorem:3.1}, we first construct approximate solutions of \eqref{eq:3.1}, 
and prove the following lemma. 
\begin{lemma}
\label{Lemma:3.1}
Assume the same conditions as in Theorem~{\rm\ref{Theorem:3.1}}. 
Then there exists a solution of \eqref{eq:3.1} such that  
\begin{gather}
\label{eq:3.6}
\sup_{0<t\le T}t^{\frac{N}{2}(1-\frac{1}{q})}\left[\|u(t)\|_q+t^{\frac{1}{2}}\|(\nabla_x u)(t)\|_q\right]<\infty,\\
\label{eq:3.7}
\sup_{0<t\le T}\left(|||u(t)|||_l+t^{\frac{1}{2}}|||(\nabla _xu)(t)|||_l\right)<\infty,
\end{gather}
for any $T>0$, $q\in[1,\infty]$, and $l\in[0,K]$. 
\end{lemma}
{\bf Proof.}
Let $q\in[1,\infty]$ and $\varphi\in L^1({\bf R}^N)$. 
Put 
\begin{equation}
\label{eq:3.8}
u_1(x,t):=(e^{t\Delta}\varphi)(x),\qquad
u_{n+1}(x,t):=(e^{t\Delta}\varphi)(x)+\int_0^t e^{(t-s)\Delta}f_n(s)ds,
\end{equation}
for $(x,t)\in{\bf R}^N\times(0,\infty)$, 
where $n=1,2,\dots$ and $f_n(y,s):=f(y,s,u_n(y,s),(\nabla u_n)(y,s))$. 
Let $c_0$ and $ c_1$ be the constants given in $(G1)$ and 
put $C:=c_0+c_1+2^{\frac{N+1}{2}}c_0c_1$.  
By $(G1)$ we have 
\begin{eqnarray}
\label{eq:3.9}
 & & \sup_{0<t<\infty}t^{\frac{N}{2}(1-\frac{1}{q})}[\|u_1(t)\|_q+t^{\frac{1}{2}}\|\nabla u_1(t)\|_q]\\
 & & = \sup_{0<t<\infty}t^{\frac{N}{2}(1-\frac{1}{q})}
 [\|e^{t\Delta}\varphi\|_q+t^{\frac{1}{2}}\|\nabla e^{t\Delta}\varphi\|_q]
\le (c_0+c_1)\|\varphi\|_1\le C\|\varphi\|_1.\nonumber
\end{eqnarray}
This together with \eqref{eq:3.2} implies that 
\begin{equation}
\label{eq:3.10}
\sup_{0<t\le T}t^{\frac{N}{2}(1-\frac{1}{q})+\frac{1}{2}}\|f_1(t)\|_q
\le CC_1(1+T)^{\frac{1}{2}}\|\varphi\|_1,\quad T>0,
\end{equation}
for some constant $C_1$. 
By $(G1)$ and \eqref{eq:3.10} we have 
\begin{eqnarray}
\label{eq:3.11}
 & & \left\|\int_0^t e^{(t-s)\Delta}f_1(s)ds\right\|_q
\le\int_0^{t/2}\|e^{(t-s)\Delta}f_1(s)\|_qds+\int_{t/2}^t\|e^{(t-s)\Delta}f_1(s)\|_qds\\
 & & \qquad\qquad
\le c_0\int_0^{t/2}(t-s)^{-\frac{N}{2}(1-\frac{1}{q})}\|f_1(s)\|_1ds
+\int_{t/2}^t\|f_1(s)\|_qds\nonumber\\
 & & \qquad\qquad
 \le CC_2(1+T)^{\frac{1}{2}}t^{-\frac{N}{2}(1-\frac{1}{q})+\frac{1}{2}}\|\varphi\|_1\nonumber
\end{eqnarray}
for all $t\in(0,T)$ and $T>0$, where $C_2$ is a constant.
Then, by $(G1)$, \eqref{eq:3.8}, and \eqref{eq:3.11} 
we have 
\begin{equation}
\label{eq:3.12}
\sup_{0<t\le T}t^{\frac{N}{2}(1-\frac{1}{q})}\|u_2(t)\|_q\le 
c_0\|\varphi\|_1+
CC_2(1+T)^{\frac{1}{2}}T^{\frac{1}{2}}\|\varphi\|_1,\quad T>0.
\end{equation}
Furthermore, since 
\begin{equation}
\label{eq:3.13}
u_2(x,t)=[e^{(t/2)\Delta}u_2(t/2)](x)+\int_{t/2}^te^{(t-s)\Delta}f_1(s)ds,
\quad (x,t)\in{\bf R}^N\times(0,\infty),
\end{equation}
applying \eqref{eq:3.10} and \eqref{eq:3.12} to \eqref{eq:3.13}, 
by $(G1)$ we obtain 
\begin{eqnarray}
\label{eq:3.14}
 & & \|\nabla u_2(t)\|_q
 \le\|\nabla e^{(t/2)\Delta}u_2(t/2)\|_q
 +\int_{t/2}^t\|\nabla e^{(t-s)\Delta}f_1(s)\|_qds\\
 & & \qquad\qquad\,\,\,\,
 \le c_1(t/2)^{-\frac{1}{2}}\|u_2(t/2)\|_q
 +c_1\int_{t/2}^t(t-s)^{-\frac{1}{2}}\|f_1(s)\|_qds\nonumber\\
 & & \qquad\qquad\,\,\,\,
 \le c_0c_1(t/2)^{-\frac{N}{2}(1-\frac{1}{q})-\frac{1}{2}}\|\varphi\|_1
+CC_3(1+T)^{\frac{1}{2}}T^{\frac{1}{2}}t^{-\frac{N}{2}(1-\frac{1}{q})-\frac{1}{2}}\|\varphi\|_1
\nonumber
\end{eqnarray}
for all $t\in(0,T)$ and $T>0$, where $C_3$ is a constant. 
Therefore, by \eqref{eq:3.12} and \eqref{eq:3.14} we have 
\begin{eqnarray}
\label{eq:3.15}
  & & \sup_{0<t\le T}t^{\frac{N}{2}(1-\frac{1}{q})}[\|u_2(t)\|_q+t^{\frac{1}{2}}\|\nabla u_2(t)\|_q]\\
 & & \le C\|\varphi\|_1+C(C_2+C_3)(1+T)^{\frac{1}{2}}T^{\frac{1}{2}}\|\varphi\|_1
 \le C\|\varphi\|_1+CC_T\|\varphi\|_1,\quad T>0,\nonumber
\end{eqnarray}
where $C_T:=(C_2+C_3)T^{\frac{1}{2}}(1+T)^{\frac{1}{2}}$. 
Furthermore we apply the same argument as in \eqref{eq:3.15} to obtain 
\begin{eqnarray*}
\sup_{0<t\le T}t^{\frac{N}{2}(1-\frac{1}{q})}[\|u_3(t)\|_q+t^{\frac{1}{2}}\|\nabla u_3(t)\|_q]
 & \!\!\!\le\!\!\! & C\|\varphi\|_1+CC_T(1+C_T)\|\varphi\|_1\\
 & \!\!\!\le\!\!\! & C(1+C_T+C_T^2)\|\varphi\|_1,\quad T>0. 
\end{eqnarray*}
Repeating the argument above, for any $n=1,2,\dots$, 
we have 
\begin{equation}
\label{eq:3.16}
\sup_{0<t\le T}t^{\frac{N}{2}(1-\frac{1}{q})}[\|u_n(t)\|_q+t^{\frac{1}{2}}\|(\nabla_x u_n)(t)\|_q]
\le C(1+C_T+\cdots+C_T^{n-1})\|\varphi\|_1
\end{equation}
and 
\begin{equation}
\label{eq:3.17}
u_{n+1}(x,t)=[e^{(t-T)\Delta}u_{n+1}(T)](x)+\int_T^te^{(t-s)\Delta}f_n(s)ds
\end{equation}
for all $(x,t)\in{\bf R}^N\times(T,\infty)$ and all $T>0$.  

Let $T_1$ be a positive constant such that $C_{T_1}\le 2^{-1}$. 
By \eqref{eq:3.16} we have 
\begin{equation}
\label{eq:3.18}
\sup_{0<t\le T_1}t^{\frac{N}{2}(1-\frac{1}{q})}[\|u_n(t)\|_q+t^{\frac{1}{2}}\|(\nabla_x u_n)(t)\|_q]
\le 2C\|\varphi\|_1.
\end{equation}
Applying the same argument as in the proof of \eqref{eq:3.18} to \eqref{eq:3.17} with $T=T_1/2$, 
we have 
$$
\sup_{T_1/2<t\le 3T_1/2}(t-T_1/2)^{\frac{N}{2}(1-\frac{1}{q})}[\|u_n(t)\|_\infty+(t-T_1/2)^{\frac{1}{2}}\|(\nabla_x u_n)(t)\|_\infty]
\le 2C\|u_n(T_1/2)\|_1
$$
for $n=1,2,\dots$. 
This together with \eqref{eq:3.18} implies that 
$$
\sup_{0<t\le 3T_1/2}t^{\frac{N}{2}(1-\frac{1}{q})}[\|u_n(t)\|_q+t^{\frac{1}{2}}\|(\nabla_x u_n)(t)\|_q]
\le C_4\|\varphi\|_1
$$
for some constant $C_4$. 
Repeating this argument, for any $T>0$, 
we can find a constant $C_5$ satisfying 
\begin{equation}
\label{eq:3.19}
\sup_{0<t\le T}t^{\frac{N}{2}(1-\frac{1}{q})}[\|u_n(t)\|_q+t^{\frac{1}{2}}\|(\nabla_x u_n)(t)\|_q]
\le C_5\|\varphi\|_1,
\quad n=1,2,\dots. 
\end{equation}
This together with \eqref{eq:3.2} implies that 
\begin{equation}
\label{eq:3.20}
 \sup_{0<t\le T}t^{\frac{N}{2}(1-\frac{1}{q})+\frac{1}{2}}\|f_n(t)\|_q\le C_6,
 \quad n=1,2,\dots. 
\end{equation}
for some constant $C_6$. 

Next, by \eqref{eq:3.20}
we apply Lemma~\ref{Lemma:2.4} and $(G1)$ to \eqref{eq:3.17}, 
and we see that, for any $\nu\in(0,1)$ and $T>0$, 
there exists a constant $C_7$, independent of $n$, such that 
\begin{equation}
\label{eq:3.21}
\frac{|u_{n+1}(x,t)-u_{n+1}(y,s)|}{|x-y|^\nu+|t-s|^{\nu/2}}
+\frac{|(\nabla_x u_{n+1})(x,t)-(\nabla_x u_{n+1})(y,s)|}{|x-y|^\nu+|t-s|^{\nu/2}}\le C_7
\end{equation}
for all $(x,t)$, $(y,s)\in{\bf R}^N\times(T/2,T)$ with $(x,t)\not=(y,s)$. 
Then, by \eqref{eq:3.19} and \eqref{eq:3.21}, 
applying the Ascoli-Arzel\`a theorem and the diagonal argument to $\{u_n\}$ 
and taking a subsequence if necessary, 
we see that there exists a function $u\in C^{\nu,\nu/2}({\bf R}^N\times(0,\infty))$ 
such that $\nabla_x u\in C^{\nu,\nu/2}({\bf R}^N\times(0,\infty))$ and 
\begin{equation}
\label{eq:3.22}
\lim_{n\to\infty}u_n(x,t)=u(x,t),\qquad
\lim_{n\to\infty}(\nabla u_n)(x,t)=(\nabla_x u)(x,t)
\end{equation}
uniformly on any compact set in ${\bf R}^N\times(0,\infty)$. 
Furthermore, by \eqref{eq:3.2}, \eqref{eq:3.19}, and \eqref{eq:3.20} we have 
\begin{equation}
\label{eq:3.23}
\sup_{0<t\le T}t^{\frac{N}{2}(1-\frac{1}{q})}[\|u(t)\|_q+t^{\frac{1}{2}}\|(\nabla_x u)(t)\|_q]<\infty,
\quad
\sup_{0<t\le T}t^{\frac{N}{2}(1-\frac{1}{q})+\frac{1}{2}}\|f(t)\|_q<\infty,
\end{equation}
for any $T>0$, where  $f(x,t)=f(x,t,u,\nabla u)$. 
In addition, we have 
\begin{equation}
\label{eq:3.24}
u(x,t)=[e^{(t-T)\Delta}u(T)](x)+\int_T^te^{(t-s)\Delta}f(s)ds
\end{equation}
for all $(x,t)\in{\bf R}^N\times(T,\infty)$ and $T>0$. 
This together with \eqref{eq:3.23} implies that $u$ is a solution of \eqref{eq:3.1}. 

It remains to prove \eqref{eq:3.7}. 
Put 
$$
w_n(t)=|||u_n(t)|||_K+t^{\frac{1}{2}}|||\nabla u_n(t)|||_K. 
$$
Then, applying $(G2)$ and $(G3)$ to \eqref{eq:3.8}, we have 
\begin{equation}
\label{eq:3.25}
\sup_{0<t<1}w_1(t)\le C_1'w_1(0)=C_1'|||\varphi|||_K<\infty
\end{equation}
for some constant $C_1'$. 
Furthermore, by \eqref{eq:3.8} we have 
\begin{eqnarray}
\label{eq:3.26}
 & & w_2(t)\le
\int_{{\bf R}^N}(1+|x|)^K(|e^{t\Delta}\varphi|+t^{\frac{1}{2}}|\nabla e^{t\Delta}\varphi|)dx
\vspace{5pt}\\ 
 & & \qquad\qquad\qquad
+\int_0^t\left(\int_{{\bf R}^N}(1+|x|)^K\left|e^{(t-s)\Delta}f_1(s)\right|dx\right)ds\nonumber\\
 & & \qquad\qquad\qquad\qquad
 +t^{\frac{1}{2}}\int_0^t\left(\int_{{\bf R}^N}(1+|x|)^K\left|\nabla e^{(t-s)\Delta}f_1(s)\right|dx\right)ds
\vspace{5pt}\nonumber\\
 & & \qquad\,\,\,
=:I_1(t)+I_2(t)+I_3(t)\nonumber
\end{eqnarray}
for all $t>0$. 
Let $T_2$ be a sufficiently small constant to be chosen later such that $0<T_2<1$. 
Then, since $I_1(t)=w_1(t)$, by \eqref{eq:3.25}
we have 
\begin{equation}
\label{eq:3.27}
\sup_{0<t\le T_2}I_1(t)\le C_1'|||\varphi|||_K.
\end{equation}
On the other hand,
by $(G2)$, \eqref{eq:3.2}, and \eqref{eq:3.27} 
we have
\begin{eqnarray}
\label{eq:3.28}
 & & I_2(t)\le C_2'\int_0^t(1+(t-s)^{\frac{K}{2}})
\left[|||f_1(s)|||_K+\|f_1(s)\|_1\right]ds\\
 & & \qquad\,\,
 \le C_3'\int_0^t s^{-\frac{1}{2}}w_1(s)ds\le C_1'C_4'T^{1/2}|||\varphi|||_K
 \nonumber
\end{eqnarray}
for all $0<t\le T<1$, where $C_2'$, $C_3'$, and $C_4'$ are constants. 
Similarly, by $(G3)$, \eqref{eq:3.2}, and \eqref{eq:3.27} 
we have 
\begin{eqnarray}
\label{eq:3.29}
 & & I_3(t)\le C_5'T^{1/2}\int_0^t (t-s)^{-\frac{1}{2}}(1+(t-s)^{\frac{K}{2}})
 \left[|||f_1(s)|||_K+\|f_1(s)\|_1\right]ds\\
 & & \qquad\,\,
 \le C_6'T^{\frac{1}{2}}\int_0^t (t-s)^{-\frac{1}{2}}s^{-\frac{1}{2}}w_1(s)ds\le C_1'C_7'T^{\frac{1}{2}}|||\varphi|||_K
 \nonumber
\end{eqnarray}
for all $0<t\le T<1$, where $C_5'$, $C_6'$, and $C_7'$ are constants.
By \eqref{eq:3.26}--\eqref{eq:3.29}, 
taking a sufficiently small $T_2>0$ so that $(C_4'+C_7')T_2^{1/2}\le 2^{-1}$, 
we have 
$$
\sup_{0<t\le T_2}w_2(t)\le C_1'[1+(C_4'+C_7')T_2^{\frac{1}{2}}]|||\varphi|||_K
\le C_1'(1+2^{-1})|||\varphi|||_K.
$$
Repeating the argument above, 
we have 
\begin{equation}
\label{eq:3.30}
\sup_{0<t\le T_2}w_n(t)\le C_1'(1+2^{-1}+\cdots+2^{-(n-1)})|||\varphi|||_K\le 2C_1'|||\varphi|||_K,
\quad n=1,2,\dots.
\end{equation}
Furthermore, applying the same argument to \eqref{eq:3.24} with $T=T_2/2$, 
by \eqref{eq:3.30} we have 
$$
\sup_{T_2/2<t\le3T_2/2}
[|||u_n(t)|||_K+(t-T_2/2)^{1/2}|||\nabla u_n(t)|||_K]
\le 2C_1'|||u_n(T_2/2)|||_K\le(2C_1')^2|||\varphi|||_K
$$
for $n=1,2,\dots$. 
This together with \eqref{eq:3.30} yields 
$$
\sup_{0<t\le 3T_2/2}w_n(t)
\le\sup_{0<t\le T_2}w_n(t)+\sup_{T_2<t\le 3T_2/2}w_n(t)\le C_8'|||\varphi|||_K<\infty,
\quad n=1,2,\dots,
$$
for some constant $C_8$. 
Repeating this argument, for any $T>0$, 
we have 
$$
\sup_{n\ge 1}\sup_{0<t\le T}w_n(t)<\infty.
$$ 
This together with \eqref{eq:3.22} implies 
$$
\sup_{0<t\le T}\left[|||u(t)|||_K+t^{\frac{1}{2}}|||(\nabla u)(t)|||_K\right]<\infty
\quad\mbox{for any $T>0$}.
$$
Thus we obtain \eqref{eq:3.7}, and the proof of Lemma~\ref{Lemma:3.1} is complete.
$\Box$\vspace{5pt}

Next we prove the following lemma. 
\begin{lemma}
\label{Lemma:3.2}
Assume the same conditions as in Theorem~{\rm\ref{Theorem:3.1}}. 
Let $u$ be a solution of \eqref{eq:3.1} given in Lemma~{\rm\ref{Lemma:3.1}}. 
Then there holds 
\begin{equation}
\label{eq:3.31}
\sup_{t>T}\,(\|u(t)\|_q+t^{\frac{1}{2}}\|\nabla u(t)\|_q)<+\infty 
\end{equation}
for any $T>0$ and $q\in[1,\infty]$. 
\end{lemma}
{\bf Proof.}
We use the same notation as in the proof of Lemma~\ref{Lemma:3.1}. 
Let $q\in[1,\infty]$. 
By \eqref{eq:3.2} we have 
\begin{equation}
\label{eq:3.32}
\|f(t)\|_q\le C_1t^{-A}(\|u(t)\|_q+t^{\frac{1}{2}}\|\nabla u(t)\|_q)
\end{equation}
for all $t\ge 1$, where $C_1$ is a constant. 
Let $T_1$ be a constant to be chosen later such that $T_1>1$. 
By $(G1)$, \eqref{eq:3.24}, and \eqref{eq:3.32}
we have 
\begin{eqnarray*}
 & & \|u(t)\|_q\le \|u(T_1)\|_q+\int_{T_1}^t \|f(s)\|_q ds\\
 & & \qquad\quad\,\,
 \le\|u(T_1)\|_q+C_1\int_{T_1}^t s^{-A}\left(\|u(s)\|_q+s^{\frac{1}{2}}\|\nabla u(s)\|_q\right)ds,\quad t\ge T_1.
 \nonumber
\end{eqnarray*}
This inequality together with $A>1$ implies that 
\begin{equation}
\label{eq:3.33}
\|u(t)\|_q\le\|u(T_1)\|_q+C_2T_1^{-A+1}\sup_{T_1\le s\le t}\left(\|u(s)\|_q+s^{\frac{1}{2}}\|\nabla u(s)\|_q\right)
\end{equation}
for all $t\ge T_1$, where $C_2$ is a constant. 
On the other hand, 
since 
\begin{eqnarray*}
 & & t^{\frac{1}{2}}\int_{T_1}^t (t-s)^{-\frac{1}{2}}s^{-A}ds
=t^{\frac{1}{2}}\left[\int_{T_1}^{t/2} (t-s)^{-\frac{1}{2}}s^{-A}ds+\int_{t/2}^t (t-s)^{-\frac{1}{2}}s^{-A}ds\right]\\
 & & 
 \le t^{\frac{1}{2}}\left[
 \left(\frac{t}{2}\right)^{-\frac{1}{2}}\int_{T_1}^{t/2}s^{-A}ds+\left(\frac{t}{2}\right)^{-A}\int_{t/2}^t(t-s)^{-\frac{1}{2}}ds\right]
\preceq T_1^{-A+1}
\end{eqnarray*}
for all $t\ge 2T_1$, 
by $(G1)$, \eqref{eq:3.24}, and \eqref{eq:3.32} 
we have 
\begin{eqnarray}
\label{eq:3.34}
 & & t^{\frac{1}{2}}\|\nabla u(t)\|_q
 \le c_1t^{\frac{1}{2}}(t-T_1)^{-\frac{1}{2}}\|u(T_1)\|_q+c_1t^{\frac{1}{2}}\int_{T_1}^t(t-s)^{-\frac{1}{2}}\|f(s)\|_q ds\\
 & & \qquad
 \le C_3\|u(T_1)\|_q+C_1c_1t^{\frac{1}{2}}
 \int_{T_1}^t(t-s)^{-\frac{1}{2}}s^{-A}\left(\|u(s)\|_q+s^{\frac{1}{2}}\|\nabla u(s)\|_q\right)ds\nonumber\\
 & & \qquad
 \le C_3\|u(T_1)\|_q+C_4T_1^{-A+1}
 \sup_{T_1\le s\le t}\left(\|u(s)\|_q+s^{\frac{1}{2}}\|\nabla u(s)\|_q\right)
 \nonumber
\end{eqnarray}
for all $t\ge 2T_1$, where $C_3$ and $C_4$ are constants independent of $T_1$. 
Let $T_1$ be a sufficiently large constant such that $C_4T_1^{-A+1}\le 1/2$. 
Then inequality \eqref{eq:3.34} together with \eqref{eq:3.6} yields  
\begin{equation}
\label{eq:3.35}
\sup_{2T_1\le s\le t}s^{\frac{1}{2}}\|\nabla u(s)\|_q
\le 2C_3\|u(T_1)\|_q+
\sup_{T_1\le s\le t}\|u(s)\|_q+
\sup_{T_1\le s\le 2T_1}s^{\frac{1}{2}}\|\nabla u(s)\|_q<\infty
\end{equation}
for all $t\ge 2T_1$. 
Furthermore, combining \eqref{eq:3.33} with \eqref{eq:3.35}, we have  
\begin{eqnarray*}
 & & \sup_{2T_1\le s\le t}\|u(s)\|_q
 \le\|u(T_1)\|_q+C_2T_1^{-A+1}\sup_{T_1\le s\le t}\left(\|u(s)\|_q+s^{\frac{1}{2}}\|\nabla u(s)\|_q\right)\\
 & & \qquad\qquad
 \le\|u(T_1)\|_q+C_2T_1^{-A+1}\sup_{T_1\le s\le 2T_1}\left(\|u(s)\|_q+s^{\frac{1}{2}}\|\nabla u(s)\|_q\right)\\
  & & \qquad\qquad
 +C_2T_1^{-A+1}\left(2C_3\|u(T_1)\|_q+\sup_{T_1\le s\le 2T_1}\|u(s)\|_q
 +\sup_{T_1\le s\le 2T_1}s^{\frac{1}{2}}\|\nabla u(s)\|_q\right)\\
 & & \qquad\qquad
 +2C_2T_1^{-A+1}\sup_{2T_1\le s\le t}\|u(s)\|_q
\end{eqnarray*}
%
%
for all $t\ge 2T_1$. 
Then, taking a sufficiently large $T_1$ so that $2C_2T_1^{-A+1}\le 1/2$ if necessary, 
we can find a constant $C_5$ satisfying 
$$
\sup_{2T_1\le s<\infty}\|u(s)\|_q\le C_5\|u(T_1)\|_q
+C_5\sup_{T_1\le s\le 2T_1}\left(\|u(s)\|_q+s^{\frac{1}{2}}\|\nabla u(s)\|_q\right)<\infty.
$$
This inequality together with \eqref{eq:3.6} implies that 
\begin{equation}
\label{eq:3.36}
\sup_{s>T}\|u(s)\|_q<\infty
\end{equation}
for any $T>0$. 
Similarly, by \eqref{eq:3.6}, \eqref{eq:3.35}, and \eqref{eq:3.36}
we have 
$\displaystyle{\sup_{s>T}}\,s^{\frac{1}{2}}\|\nabla u(s)\|_q<\infty$ for any $T>0$, 
and obtain inequality \eqref{eq:3.31}. Thus Lemma~\ref{Lemma:3.2} follows. 
$\Box$\vspace{7pt}

Now we are ready to prove Theorem~\ref{Theorem:3.1}. 
\vspace{3pt}
\newline
{\bf Proof of Theorem~\ref{Theorem:3.1}.}
Let $\varphi\in L^1_K$ with $K\ge0$. 
Let $u$ be a solution of \eqref{eq:3.1} given in Lemma~\ref{Lemma:3.1}.  
We first prove \eqref{eq:3.3}. 
Let $q\in[1,\infty]$ and assume  
\begin{equation}
\label{eq:3.37}
\sup_{t>1}t^\gamma\left(\|u(t)\|_q+t^{\frac{1}{2}}\|\nabla u(t)\|_q\right)<\infty
\end{equation}
for some $\gamma\ge 0$. 
Applying $(G1)$, \eqref{eq:3.31}, \eqref{eq:3.32}, and \eqref{eq:3.37} to inequality \eqref{eq:3.24} with $T=t/2$,
we obtain
\begin{eqnarray}
\label{eq:3.38}
 & & \|u(t)\|_q \le\|e^{(t/2)\Delta}u(t/2)\|_q+\int_{t/2}^t\|f(s)\|_qds\\
 & & \qquad\quad\,\,
\preceq t^{-\frac{N}{2}(1-\frac{1}{q})}\|u(t/2)\|_1
+\int_{t/2}^t s^{-A}(\|u(s)\|_q+s^{\frac{1}{2}}\|\nabla u(s)\|_q)ds\nonumber\\
 & & \qquad\quad\,\,
 \preceq t^{-\frac{N}{2}(1-\frac{1}{q})}+t^{-\gamma-A+1}
 \nonumber
\end{eqnarray}
for all $t\ge 2$. 
Similarly we have 
\begin{eqnarray}
\label{eq:3.39}
 & & t^{\frac{1}{2}}\|\nabla u(t)\|_q \le t^{\frac{1}{2}}\|\nabla e^{(t/2)\Delta}u(t/2)\|_q
 +t^{\frac{1}{2}}\int_{t/2}^t\|\nabla e^{(t-s)\Delta}f(s)\|_qds\\
 & & \quad
\preceq t^{-\frac{N}{2}(1-\frac{1}{q})}\|u(t/2)\|_1
+t^{\frac{1}{2}}\int_{t/2}^t 
(t-s)^{-\frac{1}{2}}s^{-A}(\|u(s)\|_q+s^{\frac{1}{2}}\|\nabla u(s)\|_q)ds\nonumber\\
 & & \quad
 \preceq t^{-\frac{N}{2}(1-\frac{1}{q})}+t^{-\gamma-A+1}
 \nonumber
\end{eqnarray}
for all $t\ge 2$. 
Then, under assumption \eqref{eq:3.37}, 
by \eqref{eq:3.31}, \eqref{eq:3.38}, and \eqref{eq:3.39}
we have 
$$
\sup_{t>1}\,t^{\kappa}\left(\|u(t)\|_q +t^{\frac{1}{2}}\|\nabla u(t)\|_q\right)<\infty,
$$
where
$$
\kappa=\min\left\{\gamma+A-1,\frac{N}{2}\left(1-\frac{1}{q}\right)\right\}.
$$
Since \eqref{eq:3.37} holds with $\gamma=0$ 
by Lemma~\ref{Lemma:3.2}, 
applying the argument above several times, 
we obtain \eqref{eq:3.37} with $\gamma=(N/2)(1-1/q)$. 
This together with \eqref{eq:3.6} implies \eqref{eq:3.3}. 

Next we prove \eqref{eq:3.4}. 
For any $l\in[0,K]$, we put 
$$
\mathcal{U}_l(t)
:=\int_{{\mathbb R}^N}|x|^l\left[|u(x,t)|+t^{\frac{1}{2}}|(\nabla_x u)(x,t)|\right]dx.
$$
Let $T$ be a sufficiently large constant to be chosen later such that $T\ge 1$. 
By \eqref{eq:3.24} we have
\begin{eqnarray}
\label{eq:3.40}
 & & \mathcal{U}_l(t)\le
\int_{{\bf R}^N}|x|^l(|e^{(t-T)\Delta}u(T)|+t^{\frac{1}{2}}|\nabla e^{(t-T)\Delta}u(T)|)dx\\
 & & \qquad\qquad
+\int_T^t\left(\int_{{\bf R}^N}|x|^l
\biggr[\left|e^{(t-s)\Delta}f(s)\right|+t^{\frac{1}{2}}\left|\nabla e^{(t-s)\Delta}f(s)\right|\biggr]dx\right)ds
\nonumber\\
 & & \qquad\,\,\,\,
 =: I_1(t)+I_2(t)\nonumber
\end{eqnarray}
for all $t>T$.
By $(G2)$, $(G3)$, and Lemma~\ref{Lemma:3.1} 
we have 
\begin{eqnarray}
\label{eq:3.41}
 & & I_1(t)\preceq\left(\int_{{\bf R}^N}|x|^l|u(x,T)|dx+(t-T)^{\frac{l}{2}}\|u(T)\|_1\right)\\
 & & \qquad\quad
 +t^{\frac{1}{2}}\left((t-T)^{-\frac{1}{2}}\int_{{\bf R}^N}|x|^l|u(x,T)|dx+(t-T)^{\frac{l-1}{2}}\|u(T)\|_1\right)
\preceq t^{\frac{l}{2}}\nonumber
\end{eqnarray}
for all $t>2T$.
Similarly,
by $(G2)$, $(G3)$, \eqref{eq:3.7}, \eqref{eq:3.32}, and Lemma~\ref{Lemma:3.2} we obtain  
\begin{eqnarray}
\label{eq:3.42}
 & & I_2(t)\preceq \int_T^t\int_{{\bf R}^N}(|y|^l+(t-s)^{\frac{l}{2}})|f(y,s)|dyds\\
 & & \qquad\qquad\qquad
 +t^{\frac{1}{2}}\int_T^t
\int_{{\bf R}^N}(|y|^l(t-s)^{-\frac{1}{2}}+(t-s)^{\frac{l-1}{2}})|f(y,s)|dyds\nonumber\\
 & & 
\preceq\int_T^t
\int_{{\bf R}^N}(|y|^l+(t-s)^{\frac{l}{2}})s^{-A}
(|u(y,s)|+s^{\frac{1}{2}}|\nabla u(y,s)|)dyds
%
%
\nonumber\\
 & & \qquad
 +t^{\frac{1}{2}}\int_T^t
\int_{{\bf R}^N}(|y|^l(t-s)^{-\frac{1}{2}}+(t-s)^{\frac{l-1}{2}})s^{-A}(|u(y,s)|+s^{\frac{1}{2}}|\nabla u(y,s)|)dyds\nonumber\\
 & & \preceq\left(\sup_{T<s<t}s^{-\frac{l}{2}}\mathcal{U}_l(s)\right)
\int_T^ts^{-A+\frac{l}{2}}ds
+\int_T^ts^{-A}(t-s)^{\frac{l}{2}}ds\nonumber\\
 & & \qquad
 +t^{\frac{1}{2}}\left(\sup_{T<s<t}s^{-\frac{l}{2}}\mathcal{U}_l(s)\right)
\int_T^ts^{-A+\frac{l}{2}}(t-s)^{-\frac{1}{2}}ds+t^{\frac{1}{2}}\int_T^ts^{-A}(t-s)^{\frac{l-1}{2}}ds\nonumber\\
 & & 
\preceq T^{-A+1}t^{\frac{l}{2}}\left(\sup_{T<s<t}s^{-\frac{l}{2}}\mathcal{U}_l(s)\right)+t^{\frac{l}{2}}\nonumber
\end{eqnarray}
for all $t>2T$.
By \eqref{eq:3.40}--\eqref{eq:3.42}
we see that there exists a constant $C_1$ such that 
$$
\sup_{2T<s<t}s^{-\frac{l}{2}}\mathcal{U}_l(s)
\le C_1T^{-A+1}\sup_{T<s<t}s^{-\frac{l}{2}}\mathcal{U}_l(s)+C_1
$$
for all $t>2T\ge 2$. 
Then, taking a sufficiently large $T$ so that $C_1T^{-A+1}\le 1/2$ if necessary, 
we have 
$$
\sup_{2T<s<\infty}s^{-\frac{l}{2}}\mathcal{U}_l(s)\le
2\sup_{T<s\le 2T}s^{-\frac{l}{2}}\mathcal{U}_l(s)+2C_1.
$$ 
This together with \eqref{eq:3.7} implies \eqref{eq:3.4}.

It remains to prove \eqref{eq:3.5}. 
Let $j=0,1$. 
For any $q\in[1,\infty]$,
by \eqref{eq:3.2} and \eqref{eq:3.3}
we have
\begin{equation}
\label{eq:3.43}
\sup_{t>0}\,(1+t)^{A-\frac{1}{2}}t^{\frac{N}{2}(1-\frac{1}{q})+\frac{1}{2}}\|f(t)\|_q<\infty. 
\end{equation}
Then, by \eqref{eq:2.3} and \eqref{eq:3.43} 
we apply Lemma~\ref{Lemma:2.2}~(i) and Lemma~\ref{Lemma:2.3}~(ii) to obtain
$$
|M_0(u(t),t)-M_0(u(t_0),t_0)|=\left|\int_{t_0}^t M_0(f(s),s)ds\right|
\preceq\int_{t_0}^t(1+s)^{-A+\frac{1}{2}}s^{-\frac{1}{2}}ds
$$
for all $t\ge t_0\ge 0$. 
This together with $A>1$ implies that
there exists a constant $M$ such that 
\begin{equation}
\label{eq:3.44}
|M_0(u(t),t)-M|=O(t^{-(A-1)})
\end{equation}
as $t\to\infty$. 
Then, by \eqref{eq:2.7} and \eqref{eq:3.44} we obtain 
\begin{equation}
\label{eq:3.45}
\lim_{t\to\infty}t^{\frac{N}{2}(1-\frac{1}{q})+\frac{j}{2}}
\|\nabla^j[M_0(u(t),t)g(t)-Mg(t)]\|_q=0
\end{equation}
for any $q\in[1,\infty]$.

Let
\begin{equation}
\label{eq:3.46}
R(x,t):=u(x,t)-M_0(u(t),t)g(t)=u(x,t)-\biggr(\int_{{\bf R}^N}u(x,t)dx\biggr)g(x,t).
\end{equation}
By Lemma~\ref{Lemma:2.3} we see that 
$$
\partial_t R=\Delta R+\tilde{f}\qquad
\mbox{in}\qquad{\bf R}^N\times(0,\infty),
$$
where
\begin{equation}
\label{eq:3.47}
\tilde{f}(x,t):=[P_0(t)f(t)](x)
=f(x,t)-\biggr(\int_{{\bf R}^N}f(x,t)dx\biggr)g(x,t).
\end{equation}
This implies that 
\begin{eqnarray*}
 & & \nabla^jR(t)=\nabla^je^{t\Delta}R(0)+\nabla^j\int_0^te^{(t-s)\Delta}\tilde{f}(s)ds\\
 & & \qquad\quad\,\,
 =\nabla^je^{t\Delta}R(0)+\biggr(\int_{t/2}^t
+\int_L^{t/2}+\int_0^L\biggr)\nabla^je^{(t-s)\Delta}\tilde{f}(s)ds\\
& & \qquad\quad\,\,
=:\nabla^je^{t\Delta}R(0)+J_1(t)+J_2(t)+J_3(t)
\end{eqnarray*}
for $t\ge2L$, where $L>0$. 
Since 
$\displaystyle{\int_{{\bf R}^N}R(x,0)dx=0}$, 
by Lemma~\ref{Lemma:2.1}~(ii) and \eqref{eq:2.7} we obtain
\begin{equation}
\label{eq:3.48}
\lim_{t\to\infty}t^{\frac{N}{2}(1-\frac{1}{q})+\frac{j}{2}}\|\nabla^je^{t\Delta}R(0)\|_q
\preceq\lim_{t\to\infty}\|e^{(t/2)\Delta}R(0)\|_1=0
\end{equation}
for any $q\in[1,\infty]$.
On the other hand, 
since it follows from 
\eqref{eq:2.7}, \eqref{eq:3.43}, and \eqref{eq:3.47} that
\begin{equation}
\label{eq:3.49}
\sup_{t>0}(1+t)^{A-\frac{1}{2}}t^{\frac{N}{2}(1-\frac{1}{q})+\frac{1}{2}}\|\tilde{f}(t)\|_q<\infty, 
\end{equation}
by $(G1)$ we have
\begin{eqnarray}
\label{eq:3.50}
t^{\frac{N}{2}(1-\frac{1}{q})+\frac{j}{2}}\|J_1(t)\|_q
& \!\!\! \preceq\!\!\!& t^{\frac{N}{2}(1-\frac{1}{q})+\frac{j}{2}}\int_{t/2}^t(t-s)^{-\frac{j}{2}}\|\tilde{f}(s)\|_qds\\
& \!\!\! \preceq \!\!\!&
t^{\frac{j}{2}-A}\int_{t/2}^t (t-s)^{-\frac{j}{2}}ds
\preceq t^{-A+1}=o(1)\nonumber
\end{eqnarray}
as $t\to\infty$.
Furthermore, by $(G1)$ and \eqref{eq:3.49}
we have 
\begin{eqnarray}
\label{eq:3.51}
 & & t^{\frac{N}{2}(1-\frac{1}{q})+\frac{j}{2}}\|J_2(t)\|_q\\
 & & 
\le t^{\frac{N}{2}(1-\frac{1}{q})+\frac{j}{2}}\int_L^{t/2}
\left\|\nabla^je^{(t-s)\Delta}\tilde{f}(s)\right\|_qds
 \preceq\int_L^{t/2}\|\tilde{f}(s)\|_1ds 
\preceq \int_L^{t/2}s^{-A}ds
\preceq L^{-A+1}\nonumber
\end{eqnarray}
for all sufficiently large $t$. 
Similarly, by $(G3)$ we have 
\begin{eqnarray}
\label{eq:3.52}
 & & t^{\frac{N}{2}(1-\frac{1}{q})+\frac{j}{2}}\|J_3(t)\|_q\\
 & & \le t^{\frac{N}{2}(1-\frac{1}{q})+\frac{j}{2}}
\int_0^L\left\|\nabla^je^{\frac{(t-s)}{2}\Delta}e^{\frac{(t-s)}{2}\Delta}\tilde{f}(s)\right\|_qds
 \preceq\int_0^L\left\|e^{\frac{(t-s)}{2}\Delta}\tilde{f}(s)\right\|_1ds
\nonumber
\end{eqnarray}
for all $t>0$. 
On the other hand, by Lemma~\ref{Lemma:2.1}~(ii), $(G1)$, \eqref{eq:2.4}, and \eqref{eq:3.49} 
we have 
\begin{gather}
\label{eq:3.53}
\lim_{t\to\infty}
\left\|e^{\frac{(t-s)}{2}\Delta}\tilde{f}(s)\right\|_1=0,\\
\label{eq:3.54}
\qquad\qquad
\left\|e^{\frac{(t-s)}{2}\Delta}\tilde{f}(s)\right\|_1
\le\|\tilde{f}(s)\|_1<\infty,\quad t\ge 2L,
\end{gather}
for all $s\in(0,L)$.  
By \eqref{eq:3.53} and \eqref{eq:3.54}
we apply the Lebesgue dominated convergence theorem to \eqref{eq:3.52}, 
and obtain 
\begin{equation}
\label{eq:3.55}
\lim_{t\to\infty}t^{\frac{N}{2}(1-\frac{1}{q})+\frac{j}{2}}\|I_3(t)\|_q=0.
\end{equation}
Therefore, by \eqref{eq:3.48}--\eqref{eq:3.51} and \eqref{eq:3.55} 
we have 
$$
\limsup_{t\to\infty}t^{\frac{N}{2}(1-\frac{1}{q})+\frac{j}{2}}\|\nabla^jR(t)\|_q
\le C_2L^{-A+1}
$$
for some constant $C_2$. 
Therefore, since $L$ is arbitrary, 
by $A>1$ we have  
$$
\lim_{t\to\infty}t^{\frac{N}{2}(1-\frac{1}{q})+\frac{j}{2}}\|\nabla^jR(t)\|_q=0. 
$$
This together with \eqref{eq:3.45} and \eqref{eq:3.46} yields \eqref{eq:3.5}, 
and Theorem \ref{Theorem:3.1} follows. 
$\Box$\vspace{3pt}

By an argument similar to the proof of Theorem \ref{Theorem:3.1} and with the aid of \eqref{eq:1.6}
%
%
we can obtain the following theorem. 
\begin{Theorem}
\label{Theorem:3.2}
Consider the Cauchy problem
\begin{equation}
\label{eq:3.56}
\partial_tu=\Delta u+\nabla\cdot\mbox{\boldmath$F$}(x,t,u)
\quad \mbox{in}\quad {\bf R}^N\times(0,\infty),\quad
u(x,0)=\varphi(x)
\quad\mbox{in} \quad {\bf R}^N,
\end{equation}
where $\mbox{\boldmath$F$}\in C({\bf R}^N\times(0,\infty)\times{\bf R}:{\bf R}^N)$ 
and $\varphi\in L^1_K$ for some $K\ge0$.
Assume that there exist constants $C>0$ and $A>1$ such that
$$
|\mbox{\boldmath$F$}(x,t,p)|\le C(1+t)^{-A+1/2}|p|,\qquad (x,t,p)\in{\bf R}^N\times(0,\infty)\times{\bf R}.
$$
Then there exists a function $u\in C({\bf R}^N\times(0,\infty))$ with the following properties: 
\begin{itemize}
  \item[{\rm(i)}] 
	For any $q\in[1,\infty]$ and $l\in[0,K]$, 
	$$	
	\sup_{0<t<\infty}t^{\frac{N}{2}(1-\frac{1}{q})}\|u(t)\|_q+\sup_{0<t<\infty}(1+t)^{-\frac{l}{2}}|||u(t)|||_l<\infty;
	$$
  \item[{\rm(ii)}]    
	$u$ satisfies  
	$$
	u(x,t)=e^{t\Delta}\varphi(x)+\int_0^t\nabla\cdot e^{(t-s)\Delta}\mbox{\boldmath$F$}(\cdot,s,u(\cdot,s))ds
	$$ 
	for almost all $(x,t)\in {\bf R}^N\times(0,\infty)$;
 \item[{\rm(iii)}]  
	There holds 
	$$
	\lim_{t\to\infty}t^{\frac{N}{2}(1-\frac{1}{q})}\|u(t)-Mg(t)\|_q=0,\quad q\in[1,\infty],
	$$
	where $M=\int_{{\bf R}^N}\varphi(x)dx$. 
\end{itemize}
\end{Theorem}
\begin{remark}
\label{Remark:3.1} 
Assume $\varphi\in L^\infty({\bf R}^N)\cap L^1_K$ for some $K\ge 0$. 
Let $u$ be the solution of \eqref{eq:3.1}, given in Theorem~{\rm\ref{Theorem:3.1}}. 
Then, by an argument similar to the proof of Lemma~\ref{Lemma:3.1} we have
%
%
$$
\sup_{0<t\le T}\left[\|u(t)\|_\infty+t^{\frac{1}{2}}\|(\nabla_x u)(t)\|_\infty\right]<\infty
$$
for any $T>0$. 
This together with assertion~{\rm (i)} of Theorem~{\rm\ref{Theorem:3.1}} implies that 
$$
\sup_{0<t<\infty}(1+t)^{\frac{N}{2}(1-\frac{1}{q})}\left[\|u(t)\|_q+t^{\frac{1}{2}}\|(\nabla_x u)(t)\|_q\right]<\infty
$$
for any $q\in[1,\infty]$. 
This also holds for the solution of \eqref{eq:3.56}, given in Theorem~{\rm\ref{Theorem:3.2}}.
\end{remark}
\section{Main Theorems}
In this section we state the main results of this paper, 
and give the higher order asymptotic expansions of the solution $u$ of Cauchy problem~\eqref{eq:1.1}.
\vspace{3pt}

Let $u$ be a solution of Cauchy problem \eqref{eq:1.1} 
with $\varphi\in L_K^1$ for some $K\ge 0$. 
Assume that the solution $u$ satisfies \eqref{eq:3.3}, \eqref{eq:3.4} and condition $(C_A)$ for some $A>1$. 
Put 
$$
F(x,t):=F(x,t,u(x,t),\nabla u(x,t))
$$ 
for simplicity. 
Then, by \eqref{eq:3.4}, 
for any multi-index $\alpha$ with $|\alpha|\le [K]$, 
we can define $M_\alpha(u(t),t)$ for all $t\ge 0$ (see \eqref{eq:2.3}). 
Furthermore, by $(C_A)$, \eqref{eq:3.3}, and \eqref{eq:3.4}
we have
\begin{eqnarray}
\label{eq:4.1}
 & & \|F(t)\|_q
 \preceq(1+t)^{-A}\left[\|u(t)\|_q+(1+t)^{\frac{1}{2}}\|\nabla_xu(t)\|_q\right]
 \preceq(1+t)^{-A+\frac{1}{2}}t^{-\frac{N}{2}(1-\frac{1}{q})-\frac{1}{2}},\\
\label{eq:4.2}
 & & |||F(t)|||_l\preceq (1+t)^{-A}\left[|||u(t)|||_l+(1+t)^{\frac{1}{2}}|||\nabla_xu(t)|||_l\right]
\preceq (1+t)^{\frac{l+1}{2}-A}t^{-\frac{1}{2}},
\end{eqnarray}
for all $t>0$, where $q\in[1,\infty]$ and $l\in[0,K]$. 
Therefore, applying Lemma~\ref{Lemma:2.2}~(i) and Lemma~\ref{Lemma:2.3}~(ii), 
we obtain 
\begin{equation}
\label{eq:4.3}
|M_\alpha(u(t),t)-M_\alpha(u(t_0),t_0)|=\left|\int_{t_0}^t M_\alpha(F(s),s)ds\right|
\preceq\int_{t_0}^t(1+s)^{-A+\frac{|\alpha|+1}{2}}s^{-\frac{1}{2}}ds
\end{equation}
for all $t\ge t_0\ge 0$. 
This implies the following: 
\newline
{\rm (i)} For any multi-index $\alpha$ with $|\alpha|\le [K]$,  
if $A>1+|\alpha|/2$,  
there exists a constant $M_\alpha$ such that 
\begin{equation}
\label{eq:4.4}
|M_\alpha(u(t),t)-M_\alpha|\preceq(1+t)^{-(A-1)+|\alpha|/2}
\quad\mbox{for all}\quad t>0;
\end{equation}
{\rm (ii)} For any multi-index $\alpha$ with $|\alpha|\le [K]$,  
if $1<A\le 1+|\alpha|/2$,  then 
\begin{equation}
\label{eq:4.5}
M_\alpha(u(t),t)=
\left\{
\begin{array}{ll}
O(t^{-(A-1)+|\alpha|/2}) & \mbox{if}\quad A<1+|\alpha|/2,\vspace{3pt}\\
O(\log t) & \mbox{if}\quad A=1+|\alpha|/2,
\end{array}
\right.
\end{equation}
as $t\to\infty$. 
\vspace{3pt}

\noindent
Now, following \cite{IK}, 
we introduce the function $U_n=U_n(x,t)$ defined inductively by
\begin{equation}
\label{eq:4.6}
\begin{array}{l}
\displaystyle
U_0(x,t):=
\sum_{|\alpha|\le [K]}M_\alpha(u(t),t)g_\alpha(x,t),
\vspace{5pt}\\
\displaystyle
U_n(x,t):=
U_0(x,t)+\int_0^te^{(t-s)\Delta}P_{[K]}(s)F_{n-1}(s)ds,
\quad n=1,2,\dots,
\end{array}
\end{equation}
where $F_{n-1}(x,t)=F(x,t,U_{n-1}(x,t), (\nabla_x U_{n-1})(x,t))$. 
In particular, since 
$$
e^{(t-s)\Delta}g_\alpha(s)=g_\alpha(t)\quad\mbox{for $t>s\ge 0$},
$$ 
by \eqref{eq:2.2} and \eqref{eq:4.6}
we have
\begin{eqnarray*}
U_n(x,t) &\!\!\!=\!\!\!& \sum_{|\alpha|\le[K]}M_\alpha(u(t),t)g_\alpha(x,t)\\
&&\qquad\quad
+\int_0^te^{(t-s)\Delta}
\biggr[F_{n-1}(s)-\sum_{|\alpha|\le[K]}M_\alpha(F_{n-1}(s),s)g_\alpha(s)\biggr]ds\\
&\!\!\!=\!\!\!& 
\sum_{|\alpha|\le[K]}\biggr[M_\alpha(u(t),t)-\int_0^t M_\alpha(F_{n-1}(s),s)ds\biggr]g_\alpha(x,t)
+\int_0^te^{(t-s)\Delta}F_{n-1}(s)ds.
\end{eqnarray*}
\vspace{3pt}

Now we are ready to state the main theorems of this paper.  
\begin{theorem}
\label{Theorem:4.1}
Let $u$ be a solution of Cauchy problem \eqref{eq:1.1} 
with $\varphi\in L_K^1$ for some $K\ge 0$. 
Assume that the solution $u$ satisfies \eqref{eq:3.3}, \eqref{eq:3.4}, and condition $(C_A)$ for some $A>1$. 
Let $n=0,1,2,\dots$ and assume condition $(F_A)$ if $n\ge 1$. 
Then there holds the following:
\vspace{3pt}
\newline
{\rm (i)} The function $U_n$ defined by \eqref{eq:4.6} satisfies 
\begin{eqnarray}
\label{eq:4.7}
&&\sup_{t>0}\,t^{\frac{N}{2}(1-\frac{1}{q})}\left[\|U_n(t)\|_q
+t^{\frac{1}{2}}\|\nabla_xU_n(t)\|_q\right]<\infty,\\
\label{eq:4.8}
&&
\sup_{t>0}\,(1+t)^{-\frac{l}{2}}\left[|||U_n(t)|||_l+t^{\frac{1}{2}}|||\nabla_xU_n(t)|||_l\right]<\infty,
\end{eqnarray}
for any $q\in[1,\infty]$ and $l\in[0,K]$;
\vspace{3pt}
\newline
{\rm (ii)} For any $q\in[1,\infty]$ and $j=0,1$,  
\begin{equation}
\label{eq:4.9}
t^{\frac{N}{2}(1-\frac{1}{q})+\frac{j}{2}}\left\|\nabla^j\bigg[u(t)-U_n(t)\bigg]\right\|_q
\preceq
\left\{
\begin{array}{l}
(1+t)^{-\frac{K}{2}}+(1+t)^{-(n+1)(A-1)}\\
\qquad\qquad\qquad\,\,\mbox{if}\quad 2(n+1)(A-1)\not=K,\vspace{5pt}\\
(1+t)^{-\frac{K}{2}}\log(2+t)\\
\qquad\qquad\qquad\,\,\mbox{if} \quad 2(n+1)(A-1)=K,
\end{array}\right.
\end{equation}
for all $t>0$;
\vspace{3pt}
\newline
{\rm (iii)} If $2(n+1)(A-1)>K$, then, for any $q\in[1,\infty]$ and $j=0,1$, 
\begin{equation}
\label{eq:4.10}
t^{\frac{N}{2}(1-\frac{1}{q})+\frac{j}{2}}\left\|\nabla^j\bigg[u(t)-U_n(t)\bigg]\right\|_q=
\left\{
\begin{array}{ll}
o(t^{-\frac{K}{2}}) & \mbox{if}\quad K=[K],\vspace{3pt}\\
O(t^{-\frac{K}{2}})& \mbox{if}\quad K>[K],
\end{array}
\right.
\end{equation}
as $t\to\infty$;
\vspace{3pt}
\newline
{\rm (iv)} 
For any $l\in[0,K]$, $\sigma>0$, and $j=0,1$, 
\begin{equation}
\label{eq:4.11}
t^{\frac{j}{2}}(1+t)^{-\frac{l}{2}}\left|\left|\left|\nabla^j\bigg[u(t)-U_n(t)\bigg]\right|\right|\right|_l
\preceq
(1+t)^{-\frac{K}{2}+\sigma}+(1+t)^{-(n+1)(A-1)}
\end{equation}
for all $t>0$.
\end{theorem}
We remark that: 
\begin{itemize}
  \item $U_n$ $(n=1,2,\dots)$ gives the $([K]+2)$-th order asymptotic expansion of the solution $u$ 
  and is determined systematically by the function $U_0$;
  \item If $2(n+1)(A-1)>K$, then the decay estimate of $\|u(t)-U_n(t)\|_q$ as $t\to\infty$ in \eqref{eq:4.10} is 
  the same as in \eqref{eq:2.5};
  \item $U_0$ is represented as a linear combination of $\{g_\alpha(x,t)\}_{|\alpha|\le [K]}$, 
  and plays a role of projection of the solution onto the space spanned by $\{g_\alpha(x,t)\}_{|\alpha|\le [K]}$. 
\end{itemize}
Furthermore we remark that 
the condition $A>1$ in Theorem~{\rm\ref{Theorem:4.1}} is crucial. 
Indeed, even if conditions $(C_A)$ and $(F_A)$ hold for some $A\in(0,1]$, 
the solution of \eqref{eq:1.1} does not necessarily behave like the Gauss kernel as $t\to\infty$,  
that is, 
the conclusions of Theorem~\ref{Theorem:4.1} does not necessarily hold. 
See Remark~\ref{Remark:6.1} and \cite[Remark~1.1]{IK}.
\vspace{3pt}

Theorem~\ref{Theorem:4.1} is an extension of \cite[Theorem 3.1]{IK}, 
and is a result for general parabolic equations. 
Next, by Theorem~\ref{Theorem:4.1} 
we give other higher order asymptotic expansions of the solution of \eqref{eq:1.1}, 
which are simple modifications of the function $U_1$. 
Let $J\in\{0,\dots,[K]\}$ and put $J_A=\min\{J,2(A-1)\}$. 
Then, by \eqref{eq:4.4} we can define the function
$$
{\mathcal U}_J(x,t):=
\left\{
\begin{array}{ll}
\displaystyle
\sum_{0\le|\alpha|<J_{A}}M_{\alpha} g_{\alpha}(x,t)
&
\mbox{if} \quad J\ge1,
\vspace{5pt}\\
\displaystyle
M g(x,t)
&
\mbox{if} \quad J=0,
\end{array}
\right. 
$$
and we write 
$F({\mathcal U}_J(x,t))=F(x,t,{\mathcal U}_J(x,t),\nabla {\mathcal U}_J(x,t))$ for simplicity. 
\begin{theorem}
\label{Theorem:4.2}
Let $u$ be a solution of Cauchy problem \eqref{eq:1.1} 
with $\varphi\in L_K^1$ for some $K\ge 0$. 
Assume that the solution $u$ satisfies \eqref{eq:3.3}, \eqref{eq:3.4}, and conditions $(C_A)$ and $(F_A)$ for some $A>1$. 
Let $J\in\{0,\dots,[K]\}$ and put 
\begin{equation}
\label{eq:4.12}
\tilde{u}(x,t):=
\sum_{|\alpha|\le[K]}M_\alpha(u(t),t)g_\alpha(x,t)+\int_0^t e^{(t-s)\Delta}P_{[K]}(s)F({\mathcal U}_J(s))ds.
\end{equation}
Then, for any $q\in[1,\infty]$ and $j=0,1$, 
\begin{equation}
\label{eq:4.13}
t^{\frac{N}{2}(1-\frac{1}{q})+\frac{j}{2}}\left\|\nabla^j\left[u(t)-\tilde{u}(t)\right]\right\|_q
=\left\{
\begin{array}{ll}
O(t^{-2(A-1)}) & \mbox{if}\quad K>4(A-1),\vspace{3pt}\\
O(t^{-\frac{K}{2}}\log t)  & \mbox{if}\quad K=4(A-1),\\
O(t^{-\frac{K}{2}}) & \mbox{if}\quad K<4(A-1),\,K\not=[K],\vspace{3pt}\\
o(t^{-\frac{K}{2}}) & \mbox{if}\quad K<4(A-1),\,K=[K],\vspace{3pt}
\end{array}
\right.
\end{equation}
as $t\to\infty$. 
\end{theorem}
Furthermore, as a corollary of Theorem~\ref{Theorem:4.2}, we have:
\begin{Corollary}
\label{Corollary:4.1}
Assume the same conditions as in Theorem $\ref{Theorem:4.1}$ and $K\ge 0$. 
Put 
\begin{eqnarray}
\label{eq:4.14}
 & & \hat{u}(x,t):=
 \sum_{|\alpha|\le[K]}M_\alpha(u(t),t)g_\alpha(x,t)+\int_0^t e^{(t-s)\Delta}P_{[K]}(s)F_M(s)ds\\
 & & \qquad\quad\,
 =\left[M-\int_0^\infty\int_{{\bf R}^N}F_M(x,t)dxdt\right]g(x,t)
 +\sum_{|\alpha|\le[K]}c_\alpha(t)g_\alpha(x,t)\nonumber\\
 & &
\qquad\qquad\qquad\qquad\qquad\qquad\qquad\qquad\qquad\qquad\quad
+\int_0^t e^{(t-s)\Delta}F_M(s)ds
 \nonumber
\end{eqnarray}
where $M=M_0$, $F_M(x,t):=F(x,t,Mg(x,t),M\nabla g(x,t))$, and  
\begin{eqnarray*}
 &&
c_0(t):=\int_t^\infty\int_{{\bf R}^N}[F_M(x,s)-F(x,s)]dxds,\\
 & &
c_\alpha(t):=M_\alpha(u(t),t)-\int_0^tM_\alpha(F_M(s),s)ds
\qquad \mbox{if}\quad1\le|\alpha|\le[K].
\end{eqnarray*} 
Then  \eqref{eq:4.13} holds with $\tilde{u}$ replaced by $\hat{u}$. 
\end{Corollary}
\section{Proof of Main Theorems}
In this section we prove Theorems~\ref{Theorem:4.1},~\ref{Theorem:4.2},
and Corollary~\ref{Corollary:4.1}. 
We first prove assertions (i), (ii), and (iv) of Theorem \ref{Theorem:4.1}. 
\vspace{3pt}
\newline
{\bf Proof of assertions (i), (ii), and (iv).} 
By \eqref{eq:3.4} we apply Lemma~\ref{Lemma:2.2}~(i) with $\beta=\gamma=0$ 
to the function $U_0$ (see \eqref{eq:4.6}), and obtain 
$$
|\nabla^jU_0(x,t)|\le\sum_{|\alpha|\le[K]}|M_\alpha(u(t),t)||\nabla^jg_\alpha(x,t)|
\preceq \sum_{|\alpha|\le[K]}(1+t)^{\frac{|\alpha|}{2}}|\nabla^jg_\alpha(x,t)|
$$
for all $(x,t)\in{\bf R}^N\times(0,\infty)$ and $j=0,1$. 
This inequality together with \eqref{eq:2.7} implies \eqref{eq:4.7} and \eqref{eq:4.8} for the case $n=0$, 
and assertion (i) follows for the case $n=0$. 

Let $n=-1,0,1,2,\dots$ and $j=0,1$. 
We assume, without loss of generality, that $\sigma\in(0,A-1)$. 
Put 
$$
\sigma_n=
\left\{
\begin{array}{ll}
\sigma & \mbox{if}\quad 2n(A-1)\ge K,\\
(K/2)-n(A-1) & \mbox{if}\quad 2n(A-1)< K,
\end{array}
\right.
\qquad
\gamma_n=A+\frac{K}{2}-\sigma_n.
$$
Let $U_{-1}\equiv 0$ and $F_{-1}\equiv 0$ in ${\bf R}^N\times(0,\infty)$. 
Then \eqref{eq:4.6} holds for $n=0,1,2,\dots$. Furthermore, 
since the solution $u$ satisfies \eqref{eq:3.3}--\eqref{eq:3.5}, 
assertions (i), (ii), and (iv) hold with $n=-1$ and $\sigma=\sigma_0$. 

We prove assertions (i), (ii), and (iv) under condition $(F_A)$. 
Assume that there exists a number $n_*\in\{-1,0,1,2,\cdots\}$ such that 
assertions (i), (ii), and (iv) hold with $n=n_*$ and $\sigma=\sigma_{n_*+1}$. 
We first prove assertion~(i) for $n=n_*+1$. 
Since $U_{n_*}\in{\cal S}$ and $0\in{\cal S}$, 
by $(F_A)$ we have 
\begin{eqnarray*}
 & & |F_{n_*}(x,t)|=|F(x,t,U_{n_*},\nabla U_{n_*})-F(x,t,0,0)|\vspace{2pt}\\
 & & \qquad\qquad\,\,\,
 \preceq(1+t)^{-A}(|U_{n_*}(x,t)|+(1+t)^{1/2}|\nabla U_{n_*}(x,t)|)
\end{eqnarray*}
for all $(x,t)\in{\bf R}^N\times(0,\infty)$. 
Then, since assertion (i) holds with $n=n_*$, 
we obtain 
$$
\sup_{t>0}\,
(1+t)^{A-\frac{1}{2}}t^{\frac{1}{2}}\left[t^{\frac{N}{2}(1-\frac{1}{q})}
\|F_{n_*}(t)\|_q+(1+t)^{-\frac{l}{2}}|||F_{n_*}(t)|||_l\right]<\infty
$$
for any $q\in[1,\infty]$ and $l\in[0,K]$. 
This together with Lemma~\ref{Lemma:2.2}~(i) implies that 
\begin{equation}
\label{eq:5.1}
\sup_{t>0}\,
(1+t)^{A-\frac{1}{2}}t^{\frac{1}{2}}\left[t^{\frac{N}{2}(1-\frac{1}{q})}
\|P_{[K]}(t)F_{n_*}(t)\|_q+(1+t)^{-\frac{l}{2}}|||P_{[K]}(t)F_{n_*}(t)|||_l\right]<\infty
\end{equation}
for any $q\in[1,\infty]$ and $l\in[0,K]$. 
Therefore, since $A>1$, by $(G1)$, \eqref{eq:4.6}, \eqref{eq:4.7} with $n=0$, and \eqref{eq:5.1}
we have  
\begin{eqnarray*}
 & & \|\nabla^jU_{n_*+1}(t)\|_q
 \le \|\nabla^jU_0(t)\|_q+\biggr\|\nabla^j\int_0^t e^{(t-s)\Delta}P_{[K]}(s)F_{n_*}(s)ds\biggr\|_q\\
 && 
 \preceq t^{-\frac{N}{2}(1-\frac{1}{q})-\frac{j}{2}}
 +\int_0^{t/2}(t-s)^{-\frac{N}{2}(1-\frac{1}{q})-\frac{j}{2}}\|P_{[K]}(s)F_{n_*}(s)\|_1ds\\
 & & \qquad\qquad\qquad\qquad\qquad\qquad
 +\int_{t/2}^t(t-s)^{-\frac{j}{2}}\|P_{[K]}(s)F_{n_*}(s)\|_qds\\
 & & \preceq t^{-\frac{N}{2}(1-\frac{1}{q})-\frac{j}{2}}
 +t^{-\frac{N}{2}(1-\frac{1}{q})-\frac{j}{2}}\int_0^{t/2}(1+s)^{-A+\frac{1}{2}}s^{-\frac{1}{2}}ds\\
 & & \qquad\qquad\qquad\qquad\qquad\qquad
 +t^{-\frac{N}{2}(1-\frac{1}{q})-\frac{1}{2}}(1+t)^{-A+\frac{1}{2}}\int_{t/2}^t(t-s)^{-\frac{j}{2}}ds\\
 & & \preceq t^{-\frac{N}{2}(1-\frac{1}{q})-\frac{j}{2}}
\end{eqnarray*}
for all $t>0$.  Furthermore, by $(G2)$, $(G3)$, \eqref{eq:4.6}, \eqref{eq:4.8} with $n=0$, and \eqref{eq:5.1} we have 
\begin{eqnarray}
\label{eq:5.2}
 & & |||\nabla^jU_{n_*+1}(t)|||_l
 \le|||\nabla^jU_0(t)|||_l+
 \biggr|\biggr|\biggr|\int_0^t \nabla^je^{(t-s)\Delta}P_{[K]}(s)F_{n_*}(s)ds \biggr|\biggr|\biggr|_l\\
 & & \preceq t^{-\frac{j}{2}}(1+t)^{\frac{l}{2}}
  +\int_0^t(t-s)^{\frac{j}{2}}|||P_{[K]}(s)F_{n_*}(s)|||_lds\nonumber\\
 & &  \qquad\qquad\qquad\qquad
 +\int_0^t(t-s)^{-\frac{j}{2}}(1+(t-s)^{\frac{l}{2}})\|P_{[K]}(s)F_{n_*}(s)\|_1ds\nonumber\\
 & & \preceq
 t^{-\frac{j}{2}}(1+t)^{\frac{l}{2}}
 +\biggr(\int_0^{t/2}+\int_{t/2}^t\biggr)(t-s)^{-\frac{j}{2}}(1+s)^{-A+\frac{l+1}{2}}s^{-\frac{1}{2}}ds\nonumber\\
 & & \qquad\qquad\qquad\,
 +\biggr(\int_0^{t/2}+\int_{t/2}^t\biggr)(t-s)^{-\frac{j}{2}}(1+(t-s)^{\frac{l}{2}})(1+s)^{-A+\frac{1}{2}}s^{-\frac{1}{2}}ds\nonumber\\
 & & \preceq t^{-\frac{j}{2}}(1+t)^{\frac{l}{2}}\nonumber
\end{eqnarray}
for all $t>0$. 
These imply that assertion (i) holds with $n=n_*+1$. 
On the other hand, 
due to $u\in{\cal S}$, by $(F_A)$ 
we have 
\begin{eqnarray}
\label{eq:5.3}
 & & |F_{n_*}(x,t)-F(x,t)|\\
 & & \preceq(1+t)^{-A}(|u(x,t)-U_{n_*}(x,t)|+(1+t)^{1/2}|\nabla u(x,t)-\nabla U_{n_*}(x,t)|)
\nonumber
\end{eqnarray}
for all $(x,t)\in{\bf R}^N\times(0,\infty)$. 
Then, since assertions (ii) and (iv) hold with $n=n_*$ and $\sigma=\sigma_{n_*+1}$, 
by \eqref{eq:5.3} we obtain 
\begin{eqnarray}
\label{eq:5.4}
& & \sup_{t>0}\,
t^{\frac{N}{2}(1-\frac{1}{q})+(\gamma_{n_*+1}-\frac{1}{2})+\frac{1}{2}}\|F(t)-F_{n_*}(t)\|_q\\
& & \qquad\qquad\quad
+\sup_{t>0}\,
(1+t)^{-\frac{l}{2}+(\gamma_{n_*+1}-\frac{1}{2})}t^{\frac{1}{2}}|||F(t)-F_{n_*}(t)|||_l<\infty\nonumber
\end{eqnarray}
for any $q\in[1,\infty]$ and $l\in[0,K]$. 
This together with Lemma~\ref{Lemma:2.2}~(i) implies that 
\begin{eqnarray}
\label{eq:5.5}
& & \sup_{t>0}\,
t^{\frac{N}{2}(1-\frac{1}{q})+(\gamma_{n_*+1}-\frac{1}{2})+\frac{1}{2}}\|P_{[K]}(t)[F(t)-F_{n_*}(t)]\|_q\\
& & \qquad\qquad\quad
+\sup_{t>0}\,
(1+t)^{-\frac{l}{2}+(\gamma_{n_*+1}-\frac{1}{2})}t^{\frac{1}{2}}|||P_{[K]}(t)[F(t)-F_{n_*}(t)]|||_l<\infty\nonumber
\end{eqnarray}
for any $q\in[1,\infty]$ and $l\in[0,K]$. 

Next we prove that assertions~(ii) and (iv) hold with $n=n_*+1$ and $\sigma=\sigma_{n_*+2}$. 
Recall that the solution $u$ satisfies \eqref{eq:3.3} and \eqref{eq:3.4}. 
Then, due to assertion~(i) with $n=n_*+1$,  
it suffices to prove that \eqref{eq:4.9} and \eqref{eq:4.11} hold with $n=n_*+1$ and $\sigma=\sigma_{n_*+2}$ 
for all sufficiently large~$t$. 
Put $z(t):=u(t)-U_{{n_*}+1}(t)$. 
Then, by \eqref{eq:2.2} and \eqref{eq:4.6} 
we have 
\begin{equation}
\label{eq:5.6}
z(x,t)=P_{[K]}(t)u(t)-\int_0^t e^{(t-s)\Delta}P_{[K]}(s)F_{n_*}(s)ds.
\end{equation}
Then, by Lemma~\ref{Lemma:2.3}~(i) 
we obtain 
$$
\partial_t z=\Delta z+P_{[K]}(t)[F(t)-F_{n_*}(t)]\qquad\mbox{in}\quad{\bf R}^N\times(0,\infty). 
$$
This implies that 
\begin{equation}
\label{eq:5.7}
z(t)=e^{(t-t_0)\Delta}z(t_0)+\int_{t_0}^t e^{(t-s)\Delta}
P_{[K]}(s)[F(s)-F_{n_*}(s)]ds,
\quad t\ge t_0\ge 0. 
\end{equation}
Let $q\in[1,\infty]$.  By $(G1)$ we have 
\begin{equation}
\label{eq:5.8}
t^{\frac{N}{2}(1-\frac{1}{q})+\frac{j}{2}}\|\nabla^j e^{t\Delta}z(0)\|_q
=t^{\frac{N}{2}(1-\frac{1}{q})+\frac{j}{2}}\|\nabla^j e^{(t/2)\Delta}e^{(t/2)\Delta} z(0)\|_q
\preceq\|e^{(t/2)\Delta} z(0)\|_1
\end{equation}
for all $t>0$. Furthermore, 
it follows from \eqref{eq:2.4} that 
$$
\int_{{\bf R}^N}x^\alpha z(x,0)dx=\int_{{\bf R}^N}x^\alpha P_{[K]}(0)u(0) dx=0,\qquad
|\alpha|\le K, 
$$
hence,
we apply \eqref{eq:5.8} and Lemma~\ref{Lemma:2.1}~(ii) to obtain
%
%
\begin{equation}
\label{eq:5.9}
t^{\frac{N}{2}(1-\frac{1}{q})+\frac{j}{2}}\|\nabla^j e^{t\Delta}z(0)\|_q
\preceq t^{-\frac{K}{2}}
\end{equation}
for all $t>0$. 
On the other hand, 
applying Lemma~\ref{Lemma:2.2}~(ii) with $\gamma'=\gamma_{n_*+1}-1/2$ and $\beta'=1/2$ 
with the aid of \eqref{eq:5.4}, we obtain  
\begin{eqnarray}
\label{eq:5.10}
 & & t^{\frac{N}{2}(1-\frac{1}{q})+\frac{j}{2}}
\left\|\nabla^j\int_0^te^{(t-s)\Delta}P_{[K]}(s)[F(s)-F_{n_*}(s)]ds\right\|_q\\
 & & \preceq
t^{-\frac{K}{2}}\int_0^t(1+s)^{\frac{K}{2}-\gamma_{n_*+1}+\frac{1}{2}}s^{-\frac{1}{2}}ds
=t^{-\frac{K}{2}}\int_0^t(1+s)^{-A+\sigma_{n_*+1}+\frac{1}{2}}s^{-\frac{1}{2}}ds\nonumber\\
 & & \preceq t^{-\frac{K}{2}}+t^{-\frac{K}{2}}\int_1^ts^{-A+\sigma_{n_*+1}}ds
 =
 \left\{
 \begin{array}{ll}
 \displaystyle{t^{-\frac{K}{2}}} & \mbox{if}\quad 2(n_*+2)(A-1)>K,\vspace{3pt}\\
 \displaystyle{t^{-\frac{K}{2}}\log t} & \mbox{if}\quad 2(n_*+2)(A-1)=K,\vspace{3pt}\\
 \displaystyle{t^{-(n_*+2)(A-1)}} & \mbox{if}\quad 2(n_*+2)(A-1)<K,
 \end{array}
 \right.\nonumber
\end{eqnarray}
for all sufficiently large $t$. 
Therefore we apply \eqref{eq:5.9} and \eqref{eq:5.10} to \eqref{eq:5.7} with $t_0=0$, 
and obtain inequality \eqref{eq:4.9} with $n=n_*+1$ for any sufficiently large $t$. 
Thus assertion~(ii) holds with $n=n_*+1$. 

On the other hand, 
for any $l\in[0,K]$, we have 
\begin{eqnarray*}
 & & (1+t)^{-\frac{l}{2}}|||\nabla^jz(t)|||_l
=\int_{{\bf R}^N}\biggr(\frac{1+|x|}{(1+t)^{1/2}}\biggr)^l|\nabla^jz(t)|dx\\
 & & 
 \preceq\int_{{\bf R}^N}\biggr[1+\biggr(\frac{1+|x|}{(1+t)^{1/2}}\biggr)^K\biggr]
 |\nabla^jz(t)|dx
 =\|\nabla^jz(t)\|_1+(1+t)^{-\frac{K}{2}}|||\nabla^jz(t)|||_K
\end{eqnarray*}
for all $t>0$. 
Then, by \eqref{eq:4.9} with $q=1$ and $n=n_*+1$ 
we see that, if there holds \eqref{eq:4.11} with $l=K$,
then we have \eqref{eq:4.11} for $l\in[0,K]$. 
Thus it suffices to prove \eqref{eq:4.11} with $l=K$, $n=n_*+1$, and $\sigma=\sigma_{n_*+2}$.
Put $Z_j(t)=|||\nabla^jz(t)|||_K$. 
By \eqref{eq:5.7} we have  
\begin{equation}
\label{eq:5.11}
Z_j(2t)\le|||\nabla^je^{t\Delta}z(t)|||_K
+\int_t^{2t}|||\nabla^je^{(2t-s)\Delta}P_{[K]}(s)[F(s)-F_{n_*}(s)]|||_Kds
\end{equation}
for all $t>0$. 
Let $\delta>0$. 
Then, by $(G2)$, $(G3)$, and \eqref{eq:4.9} with $n=n_*+1$ we have 
\begin{gather}
\label{eq:5.12}
|||e^{t\Delta}z(t)|||_K
\le (1+\delta)|||z(t)|||_K
+C_2(1+t^{\frac{K}{2}})\|z(t)\|_1
\le (1+\delta)Z_0(t)+C_3t^{\sigma_{n_*+2}},\\
\label{eq:5.13} 
t^{\frac{1}{2}}|||\nabla e^{t\Delta}z(t)|||_K
\preceq |||z(t)|||_K+(1+t^{\frac{K}{2}})\|z(t)\|_1
\preceq Z_0(t)+t^{\sigma_{n_*+2}}, 
\end{gather}
for all $t\ge 1/2$, where $C_2$ and $C_3$ constants. 
Furthermore, by $(G2)$, $(G3)$, and \eqref{eq:5.5} we have 
\begin{eqnarray}
\label{eq:5.14}
 & & \int_t^{2t}|||\nabla^je^{(2t-s)\Delta}P_{[K]}(s)[F(s)-F_{n_*}(s)]|||_K ds\\
 & & \preceq\int_t^{2t}(2t-s)^{-\frac{j}{2}}|||P_{[K]}(s)[F(s)-F_{n_*}(s)]|||_Kds\nonumber\\
 & & \qquad\qquad\qquad
 +\int_t^{2t}(2t-s)^{-\frac{j}{2}}
 \biggr[1+(2t-s)^{\frac{K}{2}}\biggr]\|P_{[K]}(s)[F(s)-F_{n_*}(s)]\|_1ds\nonumber\\
 & & \preceq\int_t^{2t}(2t-s)^{-\frac{j}{2}}(1+s)^{\frac{K}{2}-\gamma_{n_*+1}}ds\nonumber\\
 & & \qquad\qquad\qquad
 +\int_t^{2t}(2t-s)^{-\frac{j}{2}}
 \biggr[1+(2t-s)^{\frac{K}{2}}\biggr](1+s)^{-\gamma_{n_*+1}}ds\nonumber\\
 & & \preceq t^{-\frac{j}{2}+\frac{K}{2}-\gamma_{n_*+1}+1}
=t^{-\frac{j}{2}-(A-1)+\sigma_{n_*+1}}
 \preceq t^{-\frac{j}{2}+\sigma_{n_*+2}}\nonumber
 \nonumber
\end{eqnarray}
for all $t\ge 1/2$. 
Therefore, by \eqref{eq:5.11}, \eqref{eq:5.12}, and \eqref{eq:5.14} 
we can find a constant $C_4$ satisfying 
\begin{equation}
\label{eq:5.15}
Z_0(2t)\le(1+\delta)Z_0(t)+C_4t^{\sigma_{n_*+2}},\qquad t\ge 1/2.
\end{equation}
Furthermore, 
since it follows from \eqref{eq:3.4} and \eqref{eq:4.8} with $n=n_*+1$ 
that $\sup_{0<t<1}Z_0(t)<\infty$,  
we apply Lemma~\ref{Lemma:2.5} to inequality \eqref{eq:5.15}, 
and obtain 
\begin{equation}
\label{eq:5.16}
Z_0(t)\preceq t^{\sigma_{n_*+2}}
\end{equation}
for all $t\ge 1$. 
This together with \eqref{eq:5.11}, \eqref{eq:5.13}, and \eqref{eq:5.14} implies that 
\begin{equation}
\label{eq:5.17}
t^{\frac{1}{2}}Z_1(t)\preceq Z_0(t)+t^{\sigma_{n_*+2}}
\preceq t^{\sigma_{n_*+2}}
\end{equation}
for all $t\ge 1$. 
By \eqref{eq:5.16} and \eqref{eq:5.17} 
we have inequality \eqref{eq:4.11} with $n=n_*+1$, $\sigma=\sigma_{n_*+2}$ 
%
%
for any sufficiently large $t$. 
Therefore assertions~(ii) and (iv) hold with $n=n_*+1$ for all $t>0$.  
Thus, by induction we see that 
\eqref{eq:4.8}, \eqref{eq:4.9} and \eqref{eq:4.11} hold with $\sigma=\sigma_{n+1}$ for all $n=0,1,2,\dots$, 
and assertions (i), (ii), and (iv) of Theorem \ref{Theorem:4.1} follow under condition~$(F_A)$. 
Furthermore, for the case $n=0$, 
since $F_{-1}\equiv 0$, the proof of \eqref{eq:4.8}, \eqref{eq:4.9} and \eqref{eq:4.11} with $\sigma=\sigma_1$
remains true without condition $(F_A)$. 
Therefore we  obtain assertions (i), (ii), and (iv) for the case $n=0$ without condition $(F_A)$, 
and the proof of assertions (i), (ii), and (iv) is complete. 
$\Box$
\vspace{5pt}

We complete the proof of Theorem \ref{Theorem:4.1}. 
\vspace{3pt}
\newline
{\bf Proof of Theorem \ref{Theorem:4.1}.} 
It suffices to prove assertion (iii) of Theorem \ref{Theorem:4.1}. 
Since there holds \eqref{eq:4.10} for the case $K>[K]$ by Theorem~\ref{Theorem:4.1}~(ii), 
it suffices to prove \eqref{eq:4.10} for the case $K=[K]$. 
Let $K=[K]$ and assume $(F_A)$. 
Let $n\in\{0,1,2,\dots\}$ be such that 
\begin{equation}
\label{eq:5.18}
2(n+1)(A-1)>K. 
\end{equation}
Then we can take a positive constant $\sigma$ so that 
\begin{equation}
\label{eq:5.19}
\frac{K}{2}-n(A-1)<\sigma<A-1,  
\end{equation} 
and put $\epsilon:=A-1-\sigma>0$. 
By \eqref{eq:4.7} 
we see $U_n\in {\cal S}$ for $n\in\{-1,0,1,\dots\}$. 
By condition~$(F_A)$ and \eqref{eq:5.19} 
we apply Theorem \ref{Theorem:4.1}~(ii) and (iv) to obtain 
\begin{eqnarray*}
 & & t^{\frac{N}{2}(1-\frac{1}{q})}
 \|F(t)-F_{n-1}(t)\|_q+(1+t)^{-\frac{l}{2}}|||F(t)-F_{n-1}(t)|||_l\\
 & & \preceq (1+t)^{-A}\sum_{j=0,1}(1+t)^{\frac{j}{2}}\left\{t^{\frac{N}{2}(1-\frac{1}{q})}
\left\|\nabla^j\bigg[u(t)-U_{n-1}(t)\bigg]\right\|_q\right.\\
& &\qquad\qquad\qquad\quad\qquad\qquad\quad
\left.
 +(1+t)^{-\frac{l}{2}}\left|\left|\left|\nabla^j\bigg[u(t)-U_{n-1}(t)\bigg]\right|\right|\right|_l\right\}
 \nonumber\\
 & & \preceq t^{-\frac{1}{2}}(1+t)^{-A+\frac{1}{2}}[(1+t)^{-\frac{K}{2}+\sigma}+(1+t)^{-n(A-1)}]\\
  & & \preceq t^{-\frac{1}{2}}(1+t)^{-A+\frac{1}{2}-\frac{K}{2}+\sigma}
 \preceq t^{-\frac{1}{2}}(1+t)^{-\frac{K}{2}-\frac{1}{2}-\epsilon}
\end{eqnarray*}
for all $t>0$, where $q\in[1,\infty]$ and $l\in[0,K]$. 
Then, putting $\tilde{F}_{n-1}(t)=P_K(t)[F(t)-F_{n-1}(t)]$, 
by Lemma~\ref{Lemma:2.2}~(i) we have 
\begin{equation}
\label{eq:5.20}
t^{\frac{N}{2}(1-\frac{1}{q})}
 \|\tilde{F}_{n-1}(t)\|_q+(1+t)^{-\frac{l}{2}}|||\tilde{F}_{n-1}(t)|||_l\preceq t^{-\frac{1}{2}}(1+t)^{-\frac{K}{2}-\frac{1}{2}-\epsilon}
\end{equation}
for all $t>0$. 
Let $j=0,1$ and put $z_n(t)=u(t)-U_n(t)$. 
By \eqref{eq:5.7}, for any $L>0$,
we have
\begin{eqnarray}
\label{eq:5.21}
 & & \nabla^jz_n(t)
 =\nabla^je^{t\Delta}z_n(0)
+\nabla^j\int_0^t e^{(t-s)\Delta}\tilde{F}_{n-1}(s)ds\\
& & \qquad\quad\,\,\,\,
=\nabla^je^{t\Delta}z_n(0)+\left(\int_{t/2}^t +\int_L^{t/2}+\int_0^L\right)
 \nabla^je^{(t-s)\Delta}\tilde{F}_{n-1}(s)ds
\nonumber\\
& & \qquad\quad\,\,\,\,
 =:\nabla^je^{t\Delta}z_n(0)+I_1(t)+I_2(t)+I_3(t)
\nonumber
\end{eqnarray}
for $t\ge 2L$.
Since $z_n(0)=P_{[K]}(0)u(0)$, by \eqref{eq:2.4} 
we have 
$$
\int_{{\bf R}^N}x^\alpha z_n(0)dx=0,\qquad |\alpha|\le[K]=K,
$$
and by $(G1)$ and Lemma~\ref{Lemma:2.1}~(ii) 
we obtain 
\begin{equation}
\label{eq:5.22}
\lim_{t\to\infty}t^{\frac{N}{2}(1-\frac{1}{q})+\frac{K+j}{2}}\|\nabla^je^{t\Delta}z_n(0)\|_q
\preceq\lim_{t\to\infty}t^{\frac{K}{2}}\|e^{\frac{t}{2}\Delta}z_n(0)\|_1=0.
\end{equation}
On the other hand, by $(G1)$ and \eqref{eq:5.20} we have 
\begin{eqnarray}
 & & \label{eq:5.23}
t^{\frac{N}{2}(1-\frac{1}{q})+\frac{j}{2}}\|I_1(t)\|_q
\le t^{\frac{N}{2}(1-\frac{1}{q})+\frac{j}{2}}\int_{t/2}^t(t-s)^{-\frac{j}{2}}\|\tilde{F}_{n-1}(s)\|_qds\\
 & & \qquad
 \preceq t^{\frac{j}{2}-\frac{K}{2}-\epsilon-1}\int_{t/2}^t (t-s)^{-\frac{j}{2}}ds
\preceq t^{-\frac{K}{2}-\epsilon}=o(t^{-\frac{K}{2}})\nonumber
\end{eqnarray}
as $t\to\infty$. Furthermore, by Lemma~\ref{Lemma:2.1}~(ii), 
$(G1)$, \eqref{eq:2.4}, and \eqref{eq:5.20} 
we have 
\begin{eqnarray}
\label{eq:5.24}
 & & t^{\frac{N}{2}(1-\frac{1}{q})+\frac{j}{2}}\|I_2(t)\|_q
\le t^{\frac{N}{2}(1-\frac{1}{q})+\frac{j}{2}}\int_L^{t/2}
\left\|\nabla^je^{\frac{(t-s)}{2}\Delta}e^{\frac{(t-s)}{2}\Delta}\tilde{F}_{n-1}(s)\right\|_qds\\
 & & \qquad
 \preceq\int_L^{t/2}\left\|e^{\frac{(t-s)}{2}\Delta}\tilde{F}_{n-1}(s)\right\|_1ds
\preceq\int_L^{t/2}(t-s)^{-\frac{K}{2}}|||\tilde{F}_{n-1}(s)|||_Kds\nonumber\\
 & & \qquad
 \preceq t^{-\frac{K}{2}}\int_L^{t/2}s^{-1-\epsilon}ds
 \preceq t^{-\frac{K}{2}}L^{-\epsilon}\nonumber
\end{eqnarray}
for all sufficiently large $t$. 
Similarly, by $(G1)$ we have 
\begin{eqnarray}
\label{eq:5.25}
 & & t^{\frac{N}{2}(1-\frac{1}{q})+\frac{j}{2}}\|I_3(t)\|_q\\
 & & \le t^{\frac{N}{2}(1-\frac{1}{q})+\frac{j}{2}}
\int_0^L\left\|\nabla^je^{\frac{(t-s)}{2}\Delta}e^{\frac{(t-s)}{2}\Delta}\tilde{F}_{n-1}(s)\right\|_qds
 \preceq\int_0^L\left\|e^{\frac{(t-s)}{2}\Delta}\tilde{F}_{n-1}(s)\right\|_1ds
\nonumber
\end{eqnarray}
for all $t>0$. 
On the other hand, by Lemma~\ref{Lemma:2.1}~(ii), \eqref{eq:2.4}, and \eqref{eq:5.20} 
we have 
\begin{gather}
\label{eq:5.26}
\lim_{t\to\infty}
 t^{\frac{K}{2}}\left\|e^{\frac{(t-s)}{2}\Delta}\tilde{F}_{n-1}(s)\right\|_1
 =\lim_{t\to\infty}
(t-s)^{\frac{K}{2}}\left\|e^{\frac{(t-s)}{2}\Delta}\tilde{F}_{n-1}(s)\right\|_1=0,\\
\label{eq:5.27}
\left\|e^{\frac{(t-s)}{2}\Delta}\tilde{F}_{n-1}(s)\right\|_1
\preceq (t-s)^{-\frac{K}{2}}|||\tilde{F}_{n-1}(s)|||_K
\preceq t^{-\frac{K}{2}}s^{-\frac{1}{2}},\quad t\ge 2L,
\end{gather}
for all $s\in(0,L)$.  
By \eqref{eq:5.26} and \eqref{eq:5.27}
we apply the Lebesgue dominated convergence theorem to \eqref{eq:5.25}, 
and obtain 
\begin{equation}
\label{eq:5.28}
t^{\frac{N}{2}(1-\frac{1}{q})+\frac{j}{2}}\|I_3(t)\|_q=o(t^{-\frac{K}{2}})
\end{equation}
as $t\to\infty$. 
Therefore, by \eqref{eq:5.21}--\eqref{eq:5.24} and \eqref{eq:5.28} 
we see that there exists a constant $C_3$ such that  
$$
\limsup_{t\to\infty}t^{\frac{N}{2}(1-\frac{1}{q})+\frac{K+j}{2}}\|\nabla^jz_n(t)\|_q
\le C_3L^{-\epsilon}. 
$$
Then, since $L$ is arbitrary, 
we have 
$$
\lim_{t\to\infty}t^{\frac{N}{2}(1-\frac{1}{q})+\frac{K+j}{2}}\|\nabla^jz_n(t)\|_q=0. 
$$
Thus we have \eqref{eq:4.10} for the case $K=[K]$ under condition $(F_A)$. 
Furthermore, similarly as in the proof of assertions~(i), (ii), (iv), 
for the case $n=0$, we have $F_{-1}\equiv 0$, and 
the proof of \eqref{eq:4.10} with $K=[K]$ remains true without condition $(F_A)$. 
Therefore we have \eqref{eq:4.10} for the case $K=[K]$, 
and the proof of Theorem \ref{Theorem:4.1} is complete. 
$\Box$
\vspace{5pt}

Next, by arguments similar to the proof of
%
%
\cite[Theorem 5.1]{IK} and Theorem~\ref{Theorem:4.1}~(iii) 
we prove Theorem~\ref{Theorem:4.2}. 
\vspace{3pt}
\newline
{\bf Proof of Theorem \ref{Theorem:4.2}.}
Let $K\ge 0$.
By \eqref{eq:2.7} and \eqref{eq:4.4}, 
for any $q\in[1,\infty]$, $l\in[0,K]$, and $j=0,1$,
we have 
\begin{eqnarray}
\label{eq:5.29}
&&
\sup_{t>0}\left[t^{\frac{N}{2}(1-\frac{1}{q})+\gamma
+\frac{j}{2}}\|\nabla^j[U_0(t)-{\mathcal U}_J(t)]\|_q
\right.\\
&&\qquad\qquad\qquad\qquad
\left.
+(1+t)^{-\frac{l}{2}+\gamma}t^{\frac{j}{2}}|||\nabla^j[U_0(t)-{\mathcal U}_J(t)]|||_l\right]<\infty,
\nonumber
\end{eqnarray}
where $\gamma=A-1$. 
Since 
$$
|F_0(t)-F({\mathcal U}_J(t))|
\preceq(1+t)^{-A}\left\{|U_0(x,t)-{\mathcal U}_J(x,t)|
+(1+t)^{\frac{1}{2}}|\nabla\left[U_0(x,t)-{\mathcal U}_J(x,t)\right]|\right\}
$$
in ${\bf R}^N\times(0,\infty)$, 
by \eqref{eq:5.29} we have 
\begin{eqnarray}
\label{eq:5.30}
 & & \sup_{t>0}t^{\frac{N}{2}(1-\frac{1}{q})+(A+\gamma-\frac{1}{2})+\frac{1}{2}}\|F_0(t)-F({\mathcal U}_J(t))\|_q\\
 & & \qquad\qquad
 +\sup_{t>0}\,(1+t)^{-\frac{l}{2}+(A+\gamma-\frac{1}{2})}t^{\frac{1}{2}}|||F_0(t)-F({\mathcal U}_J(t))|||_l<\infty
\nonumber
\end{eqnarray}
for any $q\in[1,\infty]$ and $l\in[0,K]$. 
Then, by \eqref{eq:5.30}, 
applying Lemma~\ref{Lemma:2.2}~(ii) with $\gamma'=A+\gamma-1/2$ and $\beta'=1/2$, 
we obtain 
\begin{eqnarray}
\label{eq:5.31}
 & & t^{\frac{N}{2}(1-\frac{1}{q})+\frac{j}{2}}
\biggr\|\nabla^j\int_0^t e^{(t-s)\Delta}P_{[K]}[F_0(s)-F({\mathcal U}_J(t))]ds\biggr\|_q\\
 & &\qquad\quad
 \preceq t^{-\frac{K}{2}}\int_0^t(1+s)^{\frac{K}{2}-A-\gamma+\frac{1}{2}}s^{-\frac{1}{2}}ds\nonumber\\
 & &\qquad\quad 
 =\left\{
\begin{array}{ll}
O(t^{-\frac{K}{2}})+O(t^{-2(A-1)})
 & \mbox{if}\quad K\not=4(A-1),\\
O(t^{-\frac{K}{2}}\log t)  & \mbox{if}\quad K=4(A-1),
\end{array}
\right.
\nonumber
\end{eqnarray}
for all sufficiently large $t$. 
Furthermore, 
if $K<4(A-1)$ and $K=[K]$, then, 
by the same argument as in the proof of Theorem~\ref{Theorem:4.1}~(iii) with the aid of \eqref{eq:5.30} 
we have 
\begin{equation}
\label{eq:5.32}
t^{\frac{N}{2}(1-\frac{1}{q})+\frac{j}{2}}
\biggr\|\nabla^j\int_0^t e^{(t-s)\Delta}P_{[K]}[F_0(s)-F({\mathcal U}_J(s))]ds\biggr\|_q
=o(t^{-\frac{K}{2}})
\end{equation}
for all sufficiently large $t$. 
Therefore, since 
$$
u(t)-\tilde{u}(t)
=[u(t)-U_1(t)]+\int_0^t e^{(t-s)\Delta}P_{[K]}[F_0(s)-F({\mathcal U}_J(t))]ds, 
$$
by Theorem~\ref{Theorem:4.1}, \eqref{eq:5.31}, and \eqref{eq:5.32} 
we have 
\begin{equation}
\label{eq:5.33} 
t^{\frac{N}{2}(1-\frac{1}{q})+\frac{j}{2}}\|\nabla^j[u(t)-\tilde{u}(t)]\|_q
 =\left\{
 \begin{array}{ll}
O(t^{-2(A-1)}) & \mbox{if}\quad K>4(A-1),\vspace{3pt}\\
O(t^{-\frac{K}{2}}\log t)  & \mbox{if}\quad K=4(A-1),\\
O(t^{-\frac{K}{2}}) & \mbox{if}\quad K<4(A-1),\,K\not=[K],\vspace{3pt}\\
o(t^{-\frac{K}{2}}) & \mbox{if}\quad K<4(A-1),\,K=[K]\vspace{3pt}
  \end{array}\right.\nonumber
\end{equation}
for all sufficiently large $t$. 
Thus we obtain \eqref{eq:4.13}, and Theorem \ref{Theorem:4.2} follows. 
$\Box$\vspace{5pt}
\newline
{\bf Proof of Corollary~\ref{Corollary:4.1}}
We apply Theorem \ref{Theorem:4.2} with $J=0$. 
Then, since  
\begin{eqnarray*}
& & \tilde{u}(x,t)=
\biggr[M-\int_t^\infty\int_{{\bf R}^N}F(s)dxds\biggr]g(x,t)
+\sum_{1\le|\alpha|\le[K]}M_\alpha(u(t),t)g_\alpha(x,t)\\
 & & \qquad\quad\quad
 +\int_0^t e^{(t-s)\Delta}F_M(s)ds-g(x,t)\int_0^t\int_{{\bf R}^N}F_M(s)dxds\\
 & &\qquad\qquad\quad\quad
 -\sum_{1\le|\alpha|\le[K]}g_\alpha(x,t)\int_0^tM_\alpha(F_M(s),s)ds\\
& & \qquad\quad=
\left[M-\int_0^\infty\int_{{\bf R}^N}F_M(t)dxdt\right]g(x,t)
+\int_0^t e^{(t-s)\Delta}F_M(s)ds\\
 & &  \qquad\quad\qquad
 +\sum_{1\le|\alpha|\le[K]}\left[M_\alpha(u(t),t)-\int_0^tM_\alpha(F_M(s),s)ds
 \right]g_\alpha(x,t)\\
& & \qquad\quad\qquad\quad\quad
-\biggr[\int_t^\infty\int_{{\bf R}^N}F(s)dxds
-\int_t^\infty\int_{{\bf R}^N}F_M(s)dxds\biggr]g(x,t)=\hat{u}(x,t),\\
\end{eqnarray*}
we see that \eqref{eq:4.13} holds with $\tilde{u}$ replaced by $\hat{u}$, 
and Corollary~\ref{Corollary:4.1} follows. 
$\Box$
\section{Applications to nonlinear parabolic equations}
In this section we apply the main results of this paper, which are given in Section~4, 
to some selected nonlinear parabolic equations. 
\subsection{Convection-diffusion equation}
Consider the Cauchy problem 
for the convection-diffusion equation 
\begin{equation}
\label{eq:6.1}
\left\{
\begin{array}{ll}
\partial_t u=\Delta u+\mbox{\boldmath $a$}
\cdot\nabla(|u|^{p-1}u)\quad & \mbox{in}\quad{\bf R}^N\times(0,\infty),\vspace{3pt}\\
u(x,0)=\varphi(x)\quad & \mbox{in}\quad{\bf R}^N,
\end{array}
\right.
\end{equation}
where $N\ge 1$, $\mbox{\boldmath $a$}\in{\bf R}^N$, $p>1$, and  
$\varphi\in L^\infty({\bf R}^N)\cap L^1_K$ for some $K\ge 0$. 
Then there exists a unique bounded solution $u$ of \eqref{eq:6.1}, 
and the large time behavior of the solution $u$ has been studied in several papers 
(see for example \cite{Carpio}, \cite{DC}, \cite{DZ}, \cite{EVZ}, \cite{EZ}, \cite{K--M}, \cite{Zua}, and references therein). 
In particular, it is known that, 
if $p>1+1/N$, then the solution $u$ behaves like the Gauss kernel and \eqref{eq:1.3} holds. 

Let $p>1+1/N$. 
Then we can easily see that 
conditions $(C_A)$ and $(F_A)$ hold with 
$$
A=A_*:=\frac{N}{2}(p-1)+\frac{1}{2}>1. 
$$
Furthermore, by Theorem~\ref{Theorem:3.1} and Remark~\ref{Remark:3.1} 
we see that the unique bounded solution $u$ of \eqref{eq:6.1} satisfies \eqref{eq:3.3} and \eqref{eq:3.4}. 
These mean that all of the assertions in Section~4 hold for the solution $u$ with $A=A_*$. 
In particular, noticing that 
$$
M=\int_{{\bf R}^N}\varphi(x)dx
=\int_{{\bf R}^N}u(x,t)dx\quad\mbox{for}\quad t>0, 
$$
we have: 
\begin{Theorem}
\label{Theorem:6.1} 
Assume $p>1+1/N$ and $\varphi\in L^\infty({\bf R}^N)\cap L^1_K$ for some $K\ge 0$. 
Let $u$ be a bounded solution of {\rm\eqref{eq:6.1}} and $A=A_*$. 
Then there holds \eqref{eq:4.13} with $\tilde{u}$ replaced by 
$$
Mg(x,t)+|M|^{p-1}M\int_0^t \mbox{\boldmath $a$}\cdot\nabla e^{(t-s)\Delta}g(s)^pds
+\sum_{1\le|\alpha|\le[K]}c_\alpha(t)g_\alpha(x,t). 
$$
\end{Theorem}
Theorem~\ref{Theorem:6.1} is a direct consequence of Corollary~\ref{Corollary:4.1}. 
We remark that, for the case $K=1$, 
a result similar to Theorem~\ref{Theorem:6.1}
%
%
has been already obtained by Duro and Carpio in \cite{DC} 
(see also \cite{Zua}). 
However, as far as we know, for the case $K\not\in\{0,1\}$, 
there are no results corresponding to Theorem~\ref{Theorem:6.1} 
for the convection-diffusion equation \eqref{eq:6.1}. 
We emphasize that 
the asymptotic expansion given in Theorem~\ref{Theorem:6.1} is a simple modification of the function $U_1$,  
and Theorem~\ref{Theorem:4.1} can give the other higher order asymptotic expansions by the use of $U_n$ $(n=2,3,\dots)$. 
\begin{Remark}
\label{Remark:6.1} 
Let $1<p\le 1+1/N$ and $M\not=0$. Then, since $0<A_*\le 1$, 
we can not apply the arguments in this paper to problem~\eqref{eq:6.1}. 
On the other hand, in this case, it is known that 
the solution of \eqref{eq:6.1} does not behave like the Gauss kernel as $t\to\infty$ 
{\rm(}see for example {\rm\cite{EVZ}}, {\rm\cite{EZ}}, and {\rm\cite{K--M}}{\rm)}, 
and we can not expect that  
the assertions of Theorem~{\rm\ref{Theorem:6.1}} hold. 
\end{Remark}
The decay estimate between the solution and its asymptotic expansion 
can give the following theorem on the classification of the decay rate of $L^q$-norm of the solution $u$. 
\begin{Theorem}
\label{Theorem:6.2} 
Assume the same conditions as in Theorem~{\rm\ref{Theorem:6.1}}. 
Then the solution $u$ satisfies either 
\begin{description}
  \item[(i)] there exists an integer $d\in\{0,\dots,[K]\}$ such that, for any $q\in[1,\infty]$ and $j=0,1$, 
  $$
  \|\nabla^ju(t)\|_q\asymp t^{-\frac{N}{2}(1-\frac{1}{q})-\frac{d}{2}-\frac{j}{2}}\quad\mbox{as $t\to\infty$; or}
  $$
  \item[(ii)] for any $q\in[1,\infty]$ and $j=0,1$, 
  $$
  \lim_{t\to\infty}t^{\frac{N}{2}(1-\frac{1}{q})+\frac{[K]}{2}+\frac{j}{2}}\|\nabla^ju(t)\|_q=0. 
  $$
\end{description}
\end{Theorem}
Theorem~\ref{Theorem:6.2} is proved by the same argument as in the proof of \cite[Corollary~1.2]{IIK} 
with Theorem~\ref{Theorem:6.1}, and we leave the details of the proof to the reader.
%
%
We remark that, if $\|u(t)\|_\infty=O(t^{-(N+d)/2})$ as $t\to\infty$ for some $j\in\{1,2,\dots\}$, 
then conditions $(C_A)$ and $(F_A)$ hold with $A=A_d:=(N+d)(p-1)/2+1/2>1$ and 
all of assertions of Theorems~\ref{Theorem:4.1} and \ref{Theorem:4.2} hold with $A=A_d$. 
\subsection{Keller-Segel System}
Consider the Keller-Segel system 
of parabolic-parabolic type
\begin{eqnarray}
\label{eq:6.2}
 & & \partial_t u=\Delta u-\nabla\cdot(u\nabla v)\quad\mbox{in}\quad{\bf R}^N\times(0,\infty),\vspace{3pt}\\
\label{eq:6.3}
 & & \partial_t v=\Delta v-v+u 
\qquad\quad\mbox{in}\quad{\bf R}^N\times(0,\infty),\vspace{3pt}\\
\label{eq:6.4}
 & & u(x,0)=\varphi(x),\quad v(x,0)=\psi(x)\quad\mbox{in}\quad{{\bf R}^N},
\end{eqnarray}
where $N\ge 1$ and 
\begin{equation}
\label{eq:6.5}
\varphi,\psi, \partial_x \psi\in L^1({\bf R}^N)\cap {\mathcal B}({\bf R}^N).
\end{equation}
Here ${\mathcal B}({\bf R}^N)$ is the Banach space 
of all bounded and uniformly continuous functions on ${\bf R}^N$. 
Cauchy Problem~\eqref{eq:6.2}--\eqref{eq:6.4} is a mathematical model describing the motion of some species 
due to chemotaxis (see \cite{K-S}), 
and the asymptotics of solution $(u,v)$ of \eqref{eq:6.2}--\eqref{eq:6.4}
has been studied intensively in many papers, see for example
\cite{KS}, \cite{MR2483011}, \cite{NSY}, \cite{NY},
\cite{Yamada}, \cite{Yamada2}, and references therein. 
In particular,
it is known that, 
for any $L>0$, there exists a positive constant $\delta$ such that,  
if 
$$
\|\varphi\|_\infty\le L,\quad
\|\varphi\|_1\le\delta,\quad
\|\nabla\psi\|_1\le\delta,\quad
\|\nabla\psi\|_\infty\le\delta,
$$
then Cauchy problem \eqref{eq:6.2}--\eqref{eq:6.4} has a unique classical solution $(u,v)$ satisfying 
\begin{equation}
\label{eq:6.6}
\sup_{t>0}\,(\|u(t)\|_p+\|v(t)\|_p)<\infty\qquad \mbox{for}\quad p\in\{1,\infty\}.
\end{equation}
(See \cite[Theorem 1.2]{NSY}.) 

Let $(u,v)$ be a classical solution of \eqref{eq:6.2}--\eqref{eq:6.4} satisfying \eqref{eq:6.6}. 
Assume $\varphi\in L^1_K$ for some $K\ge 0$. 
Then we show that higher order asymptotic expansions of the solution of \eqref{eq:6.2}--\eqref{eq:6.4} 
are given as a corollary of our results. 
%
%
By \cite[Proposition~4.1]{NY} we have  
\begin{eqnarray}
\label{eq:6.7}
 & & \sup_{t>0}\,(1+t)^{\frac{N}{2}(1-\frac{1}{q})}\|u(t)\|_q+
\sup_{t\ge 1}\,t^{\frac{N}{2}(1-\frac{1}{q})+\frac{1}{2}}\|\nabla u(t)\|_q\\
 & & \qquad\qquad
 +\sup_{t>0}\,(1+t)^{\frac{N}{2}(1-\frac{1}{q})+\frac{1}{2}}\|\nabla v(t)\|_q<\infty
 \quad\mbox{for any $q\in[1,\infty]$}.
\nonumber
\end{eqnarray}
Furthermore, applying arguments similar to the proof \cite[Proposition~4.1]{NY},
%
%
we can easily obtain
\begin{equation}
\label{eq:6.8}
\sup_{t\ge 1}\,t^{\frac{N}{2}(1-\frac{1}{q})+1}\|\nabla^2 u(t)\|_q<\infty
\quad\mbox{for any $q\in[1,\infty]$}. 
\end{equation}
In addition, by \eqref{eq:6.7} 
we can apply Theorem~\ref{Theorem:3.2} to \eqref{eq:6.2}, and see that 
the solution $u$ satisfies all of the assertions of Theorem~\ref{Theorem:3.2}. 
On the other hand, since it follows from \eqref{eq:6.3} that 
\begin{equation}
\label{eq:6.9}
v(t)=e^{-t}e^{t\Delta}\psi+\int_0^t e^{-t+s}e^{(t-s)\Delta}u(s)ds,
\qquad t>0,
\end{equation}
by $(G1)$, \eqref{eq:6.7}, and \eqref{eq:6.8} we have 
\begin{equation}
\label{eq:6.10}
\sup_{t\ge 1}\,t^{\frac{N}{2}(1-\frac{1}{q})+1}\|\nabla^2 v(t)\|_q<\infty
\quad\mbox{for any $q\in[1,\infty]$}. 
\end{equation}
Therefore, putting 
\begin{equation}
\label{eq:6.11}
F(x,t,u,\nabla u):=-\nabla\cdot(u\nabla v)
=-\nabla v\cdot\nabla u-(\Delta v)u,
\end{equation}
by \eqref{eq:6.7} and \eqref{eq:6.10} 
we see that, in \eqref{eq:6.2}, 
there hold conditions $(C_A)$ and $(F_A)$ in ${\bf R}^N\times(1,\infty)$ with 
$$
A=\frac{N}{2}+1\ge\frac{3}{2}. 
$$
Furthermore, by Theorem~\ref{Theorem:3.2}~(i)
we have $u(1)\in L^1_K$. Therefore, taking the function $u(1)$ as the initial function 
of parabolic equation~\eqref{eq:6.2}, we see that all of the assertions in Section~4 hold 
with $A=N/2+1$ for the solution $u$. 
In particular, we have
\begin{Lemma} 
\label{Lemma:6.1}
Let $(u,v)$ be a global  in time solution of \eqref{eq:6.2}--\eqref{eq:6.4} satisfying \eqref{eq:6.6}. 
Assume $\varphi\in L_K^1$ for some $K\ge 0$.
Let $c_\alpha(t)$ be the functions given in Corollary~$\ref{Corollary:4.1}$.
Then there holds the following: 
\begin{itemize}
  \item[{\rm (a)}]  $c_0(t)=0$ for all $t>0$;
  \item[{\rm(b)}]  If $|\alpha|\le[K]$ and $1\le|\alpha|<N$, then 
  there exists a constant $c_\alpha$ such that 
  $$
  c_\alpha(t)=c_\alpha+O(t^{-\frac{N}{2}+\frac{|\alpha|}{2}})\quad\mbox{as}\quad t\to\infty;
  $$
  \item[{\rm(c)}] If $|\alpha|\le[K]$ and $1\le|\alpha|=N$, then 
  $c_\alpha(t)=O(\log t)$ as $t\to\infty$;
  \item[{\rm(d)}] 
  $$
  t^{\frac{N}{2}(1-\frac{1}{q})+\frac{j}{2}}\biggr\|\nabla^j\int_0^t e^{(t-s)\Delta}F_M(s)ds\biggr\|_q
  =O(t^{-\frac{N}{2}})\quad\mbox{as}\quad t\to\infty.
  $$
\end{itemize}
\end{Lemma}
{\bf Proof.} Assertion~(a) follows from \eqref{eq:6.11} and the definition of $c_0(t)$. 
Furthermore, since 
$$
\sup_{t>0}|M_\alpha(f,t)|\preceq |||f|||_{|\alpha|}\quad\mbox{for}\quad f\in L^1({\bf R}^N,(1+|x|)^{|\alpha|}dx), 
$$
by \eqref{eq:2.7}, \eqref{eq:6.7}, \eqref{eq:6.10}, and \eqref{eq:6.11} we have 
$$
|M_\alpha(F_M(t),t)|
\preceq \|\nabla v(t)\|_\infty|||\nabla g(t)|||_{|\alpha|}+\|\Delta v(t)\|_\infty|||g(t)|||_{|\alpha|}\\
\preceq t^{-\frac{N}{2}-1+\frac{|\alpha|}{2}}
$$
for all sufficiently large $t$. Then, by using \eqref{eq:4.4} and \eqref{eq:4.5} with $A=N/2+1$ 
we have assertions~(b) and (c). In addition, by $(G1)$, \eqref{eq:6.7}, \eqref{eq:6.10}, and \eqref{eq:6.11} we have 
\begin{eqnarray*}
 & & t^{\frac{N}{2}(1-\frac{1}{q})+\frac{j}{2}}\biggr\|\nabla^j\int_0^t e^{(t-s)\Delta}F_M(s)ds\biggr\|_q\\
 & & \preceq\int_0^{t/2}\|F_M(s)\|_1ds
+t^{\frac{N}{2}(1-\frac{1}{q})}\int_{t/2}^t (t-s)^{-\frac{j}{2}}\|F_M(s)\|_qds\\
 & & \preceq\int_0^{t/2}(1+s)^{-\frac{N}{2}-1}ds
 +t^{\frac{N}{2}(1-\frac{1}{q})+\frac{j}{2}}\int_{t/2}^t (t-s)^{-\frac{j}{2}}s^{-\frac{N}{2}-1-\frac{N}{2}(1-\frac{1}{q})}ds
\preceq t^{-\frac{N}{2}}
\end{eqnarray*}
for all sufficiently large $t$. This gives assertion~(d), and Lemma~\ref{Lemma:6.1} follows.
$\Box$\vspace{3pt}
\newline
Then, since 
$$
M\equiv \int_{{\bf R}^N}\varphi(x)dx=\int_{{\bf R}^N}u(x,t)dx\quad\mbox{for}\quad t>0,  
$$
by Lemma~\ref{Lemma:6.1} 
we apply Corollary~\ref{Corollary:4.1} with $N\ge K$ to obtain the following theorem.  
\begin{Theorem}
\label{Theorem:6.3} 
Let $(u,v)$ be a global  in time solution of \eqref{eq:6.2}--\eqref{eq:6.4}, satisfying \eqref{eq:6.6}. 
Let $N\ge K$ and assume $\varphi\in L_K^1$. 
Then, for any $j=0,1$, there holds the following:
\newline
{\rm (i)} If $N>K$, then  
\begin{equation}
\label{eq:6.12}
t^{\frac{N}{2}(1-\frac{1}{q})+\frac{j}{2}}
\biggr\|\nabla^j\biggr[u(t)-Mg(t)-\sum_{1\le|\alpha|\le[K]}c_\alpha g_\alpha(t)\biggr]\biggr\|_q
=\left\{
\begin{array}{ll}
o(t^{-\frac{K}{2}}) & \mbox{if}\quad K=[K],\vspace{3pt}\\
O(t^{-\frac{K}{2}})& \mbox{if}\quad K>[K],
\end{array}
\right.
\end{equation}
as $t\to\infty$;
\vspace{3pt}
\newline
{\rm (ii)} if $N=K$, then
$$
t^{\frac{N}{2}(1-\frac{1}{q})+\frac{j}{2}}
\biggr\|\nabla^j\biggr[u(t)-Mg(t)-\sum_{1\le|\alpha|\le K -1}c_\alpha g_\alpha(t)
-\sum_{|\alpha|=K}c_\alpha(t) g_\alpha(t)\biggr]\biggr\|_q=o(t^{-\frac{K}{2}})
$$
and
\begin{equation}
\label{eq:6.13}
t^{\frac{N}{2}(1-\frac{1}{q})+\frac{j}{2}}
\biggr\|\nabla^j\biggr[u(t)-Mg(t)-\sum_{1\le|\alpha|\le K -1}c_\alpha g_\alpha(t)]\biggr\|_q=O(t^{-\frac{K}{2}}\log t),
\end{equation}
as $t\to\infty$;
\vspace{3pt}
\newline
{\rm (iii)} if $N=K=1$, then
\begin{equation}
\label{eq:6.14}
t^{\frac{1}{2}(1-\frac{1}{q})+\frac{j}{2}}
\left\|\nabla^j[u(t)-Mg(t)]\right\|_q=O(t^{-\frac{1}{2}})\quad\mbox{as}\quad t\to\infty;
\end{equation}
{\rm (iv)} The same assertions as in \eqref{eq:6.12}--\eqref{eq:6.14} hold for $v$.
\end{Theorem}
{\bf Proof of Theorem \ref{Theorem:6.3}.}
Assertions~(i) and (ii) follow from Corollary~\ref{Corollary:4.1} and Lemma~\ref{Lemma:6.1}. 
Furthermore, by \eqref{eq:6.9} we see that 
\eqref{eq:6.12} and \eqref{eq:6.13} hold with $u$ replaced by $v$. 

We prove assertion~(iii). 
For this aim, by \eqref{eq:2.7} and assertion~(ii) we have only to prove 
\begin{equation}
\label{eq:6.15}
c_\alpha(t)=O(1)\quad\mbox{as}\quad t\to\infty
\end{equation}
for the case $K=N=|\alpha|=1$.  
Since $\displaystyle{\int_{\bf R}g g_xdx=0}$ and \eqref{eq:6.13} hold for $u$ and $v$, 
by \eqref{eq:2.3} and \eqref{eq:6.7} we have
\begin{eqnarray*}
 & & \left|M_\alpha(F(t),t)\right|
=\left|\int_{\bf R}x(u(x,t) v_x(x,t))_xdx\right|
=\left|\int_{\bf R} u(x,t)v_x(x,t)dx\right|\\
 & &
 =\left|\int_{\bf R}u(x,t)(v(x,t)-Mg(x,t))_xdx\right|+\left|\int_{\bf R}(Mg(x,t))_x(u(x,t)-Mg(x,t))dx\right|\\
 & &
 \le\|u(t)\|_\infty\|(v(t)-Mg(t))_x\|_1+\|(Mg(t))_x\|_\infty\|u(t)-Mg(t)\|_1
 =o(t^{-\frac{3}{2}}\log t) 
\end{eqnarray*}
as $t\to\infty$. 
Similarly we have 
$$
\left|M_\alpha(F_M(t),t)\right|
  =\left|\int_{\bf R}x(Mg(t)v_x(x,t))_xdx\right|
=\left|\int_{\bf R} Mg(x,t)v_x(x,t)dx\right|=o(t^{-\frac{3}{2}}\log t)
$$
as $t\to\infty$. These together with Lemma~\ref{Lemma:2.3}~(ii) implies \eqref{eq:6.15}, 
and assertion~(iii) follows. Then, by \eqref{eq:6.9} we see that 
\eqref{eq:6.14} holds with $u$ replaced by $v$,  
and Theorem~\ref{Theorem:6.3} follows.
$\Box$\vspace{5pt}
\begin{Remark}
\label{Remark:6.2}
{\rm (i)} Under assumption \eqref{eq:6.6}, 
Kato in {\rm\cite{MR2483011}} and Yamada in {\rm\cite{Yamada}} and {\rm\cite{Yamada2}} 
recently studied the asymptotic expansions of the solution 
of \eqref{eq:6.2}--\eqref{eq:6.4} in detail, 
and obtained some asymptotic expansions given in Theorem~{\rm\ref{Theorem:6.3}}. 
We emphasize that 
Theorem~{\rm\ref{Theorem:6.3}} is easily obtained by Corollary~{\rm\ref{Corollary:4.1}} with the aid of 
some global bounds of the solution and that 
Theorems~{\rm\ref{Theorem:4.1}} and {\rm\ref{Theorem:4.2}} can systematically give 
the other higher order asymptotic expansions of the solution and the decay estimates 
between the solution and its asymptotic expansions. 
\vspace{3pt}
\newline
{\rm (ii)} Due to the decay estimates in Theorem~{\rm\ref{Theorem:6.3}}, 
we can obtain the result similar to Theorem~{\rm\ref{Theorem:6.2}},
%
%
and by using Theorems~{\rm\ref{Theorem:4.1}} and {\rm\ref{Theorem:4.2}} 
we can also give the higher order asymptotic expansions of the solutions decaying faster than the Gauss kernel. 
\end{Remark}
\subsection{System of semilinear parabolic equations}
Our arguments in this paper are also
applicable to systems of parabolic equations 
under suitable assumptions. In this subsection we focus on the Cauchy problem 
for a system for semilinear parabolic equations, 
\begin{equation}
\label{eq:6.16}
\partial_t\mbox{\boldmath$u$}=\Delta\mbox{\boldmath$u$}+\mbox{\boldmath$F$}(\mbox{\boldmath$u$})
\quad\mbox{in}\quad{\bf R}^N\times(0,\infty),
\qquad
\mbox{\boldmath$u$}(x,0)=\Phi(x)
\quad\mbox{in}\quad{\bf R}^N,
\end{equation}
where $m=1,2,\dots$, $\mbox{\boldmath$u$}=(u_1,\cdots,u_m)$, 
$\mbox{\boldmath$F$}=(F_1(\mbox{\boldmath$u$}),\cdots F_m(\mbox{\boldmath$u$}))$, and 
$\Phi=(\varphi_1,\cdots,\varphi_m)\in (L_K^1\cap L^\infty({\bf R}^N))^m$ for some $K\ge 0$, 
and we study the asymptotics of the solution $\mbox{\boldmath$u$}$. 
Throughout this subsection we assume $\mbox{\boldmath$F$}\in C({\bf R}^N:{\bf R}^m)$ and 
that there exist constants $C>0$ and $a>1+2/N$ such that 
\begin{equation}
\label{eq:6.17}
|\mbox{\boldmath$F$}(\mbox{\boldmath$v$})|\le C|\mbox{\boldmath$v$}|^a,\qquad \mbox{\boldmath$v$}\in{\bf R}^m. 
\end{equation}

Let $\mbox{\boldmath$u$}$ be a unique global in time solution of \eqref{eq:6.16} such that 
\begin{equation}
\label{eq:6.18}
\|\mbox{\boldmath$u$}(t)\|_\infty\preceq (1+t)^{-\frac{N}{2}},\qquad t>0. 
\end{equation} 
Then, by \eqref{eq:6.17} and \eqref{eq:6.18} we have 
\begin{equation}
\label{eq:6.19}
|\mbox{\boldmath$F$}(\mbox{\boldmath$u$}(x,t))|\preceq (1+t)^{-\frac{N(a-1)}{2}}|\mbox{\boldmath$u$}(x,t)|
\end{equation} 
for all $(x,t)\in{\bf R}^N\times(0,\infty)$. 
Therefore, similarly to Section~6.1, 
we can apply the same arguments as in the previous sections  
to the solution $\mbox{\boldmath$u$}$ with $A=N(a-1)/2>1$. 
This means that all of the assertions in Section~4 hold with $A=N(a-1)/2>1$. 
In particular, we apply Corollary~\ref{Corollary:4.1} with $K\in(0,1]$ to obtain the following result.
This is an extension of \cite[Theorem~5.1]{IK}, which treats the case $m=1$.
\begin{theorem}
\label{Theorem:6.4}
Let $m\in\{1,2,\dots\}$ and $K\ge 0$. 
Assume \eqref{eq:6.17} and $\Phi=(\varphi_1,\cdots,\varphi_m)\in (L_K^1\cap L^\infty({\bf R}^N))^m$. 
Let $\mbox{\boldmath$u$}$ be a global in time solution of Cauchy problem \eqref{eq:6.16}, satisfying \eqref{eq:6.18}. 
Then there exists the limit 
$$
\mbox{\boldmath$M$}:=\lim_{t\to\infty}\int_{{\bf R}^N}\mbox{\boldmath$u$}(x,t)dx
$$
such that 
$$
\lim_{t\to\infty}t^{\frac{N}{2}(1-\frac{1}{q})}\|\mbox{\boldmath$u$}(t)-\mbox{\boldmath$M$}g(t)\|_q=0
$$
for any $q\in[1,\infty]$. Furthermore there holds the following: 
\vspace{3pt}
\newline
{\rm (i)} If $K\in(0,1]$, then 
$$
t^{\frac{N}{2}(1-\frac{1}{q})}\|\mbox{\boldmath$u$}(t)-\mbox{\boldmath$M$}g(t)\|_q
=\left\{
\begin{array}{ll}
O(t^{-\frac{K}{2}})+O(t^{-(A-1)}) & \mbox{if}\quad 2(A-1)\not=K,\vspace{5pt}\\
O(t^{-\frac{K}{2}}\log t) & \mbox{if} \quad 2(A-1)=K,
\end{array}\right.
$$
as $t\to\infty$, for any $q\in[1,\infty]$\vspace{5pt};
\newline
{\rm(ii)}
If $K\in(0,1]$, then 
$$
t^{\frac{N}{2}(1-\frac{1}{q})}\|\mbox{\boldmath$u$}(t)-\mbox{\boldmath$u$}_1(t)\|_q
=\left\{
\begin{array}{ll}
O(t^{-\frac{K}{2}})+O(t^{-2(A-1)}) & \mbox{if}\quad 4(A-1)\not=K,\vspace{5pt}\\
O(t^{-\frac{K}{2}+\sigma}) & \mbox{if} \quad 4(A-1)=K,
\end{array}\right.
$$
as $t\to\infty$, for any $q\in[1,\infty]$ and $\sigma>0$, where 
$$
\mbox{\boldmath$u$}_1(x,t)
=\left(\mbox{\boldmath$M$}-\int_0^\infty\int_{{\bf R}^N}\mbox{\boldmath$F$}(\mbox{\boldmath$M$}g(x,t))dxdt\right)g(x,t)
+\int_0^te^{(t-s)\Delta}\mbox{\boldmath$F$}(\mbox{\boldmath$M$}g(x,s))ds;
$$
{\rm(iii)}
Assume that 
$\displaystyle{\int_{{\bf R}^N}x\mbox{\boldmath$F$}(\mbox{\boldmath$M$}g(t))dx=0}$ for all $t>0$. 
Let $K>1$. 
Then 
$$
t^{\frac{N}{2}(1-\frac{1}{q})}\|\mbox{\boldmath$u$}(t)-\mbox{\boldmath$M$}g(t)\|_q
= O(t^{-\frac{1}{2}})+O(t^{-(A-1)})
$$
as $t\to\infty$ for any $q\in[1,\infty]$.
\end{theorem}
{\bf Proof of Theorem~\ref{Theorem:6.4}.}
This theorem is proved by Corollary~\ref{Corollary:4.1} with minor modifications. 
We leave the details of the proof to the reader.
%
%
(See also the proof of \cite[Proposition~5.1]{IK}.)
$\Box$
\section{Appendix}
For convenience we present the proof of Lemma~\ref{Lemma:2.4}  
by the same arguments as in Chapter 1 in \cite{F}. 
We first prove \eqref{eq:2.9} and \eqref{eq:2.10}.
\vspace{5pt}

\noindent
{\bf Proof of \eqref{eq:2.9} and \eqref{eq:2.10}.}
The $C^1$-regularity of $w$ and the representation \eqref{eq:2.9} 
are easily obtained by a argument similar to Chapter 1 of \cite{F}.
%
%
Put $C_H=\|H\|_{L^\infty(0,T:L^\infty({\bf R}^N))}$. 
Then, by \eqref{eq:2.6}, \eqref{eq:2.8}, and \eqref{eq:2.9} we see that 
there exist constants $C_1$, $C_2$, and $C_3$, independent of $C_H$ and $T$, such that 
\begin{eqnarray*}
 & & |w(x,t)|\le\int_0^t\left(\int_{{\bf R}^N}G(x-\xi,\tau)d\xi\right)\|H(\tau)\|_\infty d\tau
\le\int_0^t\|H(\tau)\|_\infty d\tau
\le C_1C_HT,\\
 & & |(\nabla_x w)(x,t)|\le\int_0^t\left(\int_{{\bf R}^N}|(\nabla_x G)(x-\xi,\tau)|d\xi\right)\|H(\tau)\|_\infty d\tau\\
 & & \qquad\qquad\qquad\!
 \le C_2\int_0^t(t-\tau)^{-\frac{1}{2}}\|H(\tau)\|_\infty d\tau\le C_3C_HT^{1/2}
\end{eqnarray*}
for all $(x,t)\in{\bf R}^N\times(0,T)$, and we obtain \eqref{eq:2.10}. 
$\Box$\vspace{5pt}

Next we prove \eqref{eq:2.11}. For this aim, 
we prove the following lemmas. Put 
$G_\alpha(x,t)=(\partial^\alpha_xG)(x,t)$.  
\begin{lemma}
\label{Lemma:7.1}
Let $0<\nu<1$ and $|\alpha|\le 1$. Then there exists a constant $C$ such that 
\begin{equation}
\label{eq:7.1}
\Pi_1(x,y:t):=
\frac{|G_\alpha(x,t)-G_\alpha(y,t)|}{|x-y|^\nu}
\le C\{h(x,t)+h(y,t)\}
\end{equation}
for all $x$, $y\in{\bf R}^N$ with $x\not=y$ and all $t>0$, where 
\begin{equation}
\label{eq:7.2}
h(x,t)=t^{-\frac{N}{2}-\frac{|\alpha|+\nu}{2}}\left[1+(t^{-\frac{1}{2}}|x|)^{-\nu}+(t^{-\frac{1}{2}}|x|)^{|\alpha|+2}\right]
e^{-\frac{|x|^2}{16t}}.
\end{equation}
\end{lemma}
{\bf Proof.}
Let $x$, $y\in{\bf R}^N$ with $x\not=y$ and $t>0$. 
If $|x-y|\ge t^{1/2}$, then, by \eqref{eq:2.6} we have 
$$
\Pi_1(x,y:t)
\le t^{-\frac{\nu}{2}}\left\{|G_\alpha(x,t)|+|G_\alpha(y,t)|\right\}
\le C_1[h(x,t)+h(y,t)]
$$
for some constant $C_1$, 
and obtain inequality \eqref{eq:7.1}. 
So it suffices to prove inequality \eqref{eq:7.1} for the case $|x-y|<t^{1/2}$. 
In this case, if $y\in B(x,|x|/2)$, 
the mean value theorem implies the existence of the point $x_*\in B(x,|x|/2)$ such that 
$$
\Pi_1(x,y:t)\le|(\nabla_x G_\alpha)(x_*,t)||x-y|^{1-\nu}\le t^{\frac{1-\nu}{2}}|(\nabla_x G_\alpha)(x_*,t)|. 
$$
Then, since $|x|/2\le|x_*|\le 3|x|/2$, 
by \eqref{eq:2.6} we have 
\begin{eqnarray}
\label{eq:7.3}
 & & \Pi_1(x,y:t) \le C_2t^{-\frac{N}{2}-\frac{|\alpha|+\nu}{2}}
 \left[1+(t^{-\frac{1}{2}}|x_*|)^{|\alpha|+1}\right]e^{-\frac{|x_*|^2}{4t}}\\
 & & \le C_3t^{-\frac{N}{2}-\frac{|\alpha|+\nu}{2}}\left[1+(t^{-\frac{1}{2}}|x|)^{|\alpha|+1}\right]e^{-\frac{|x|^2}{16t}}
\le C_4h(x,t)\quad\mbox{if}\quad y\in B(x,|x|/2),
\nonumber
\end{eqnarray}
where $C_2$, $C_3$, and $C_4$ are constants independent of $x$, $y$ and $t$. 
Similarly we have 
\begin{equation}
\label{eq:7.4}
\Pi_1(x,y:t)\le C_4h(y,t)
\quad\mbox{if}\quad x\in B(y,|y|/2).
\end{equation}
On the other hand, if $y\not\in B(x,|x|/2)$ and $x\not\in B(y,|y|/2)$, 
then we have 
$$
|x-y|\ge(1/2)\min\{|x|,|y|\},
$$ 
and obtain 
$$
\Pi_1(x,y:t) 
\le t^{-\frac{\nu}{2}}
\left[(t^{-\frac{1}{2}}|x|)^{-\nu}|G_\alpha(x,t)|
+(t^{-\frac{1}{2}}|y|)^{-\nu}|G_\alpha(y,t)|\right].
$$
This together with \eqref{eq:2.6} implies that 
\begin{equation}
\label{eq:7.5}
\Pi_1(x,y:t)\le C_5[h(x,t)+h(y,t)], 
\end{equation}
where $C_5$ is a constant independent of $x$, $y$ and $t$. 
Therefore, by \eqref{eq:7.3}--\eqref{eq:7.5} we have 
inequality \eqref{eq:7.1} for the case $|x-y|\le t^{1/2}$.
Thus Lemma~\ref{Lemma:7.1} follows. 
$\Box$
\begin{lemma}
\label{Lemma:7.2}
Let $0<\nu<1$ and $|\alpha|\le 1$. Then there exists a constant $C$ such that 
\begin{equation}
\label{eq:7.6}
\Pi_2(t,s:x):=
\frac{|G_\alpha(x,t)-G_\alpha(x,s)|}{|t-s|^{\nu/2}}
\le C\{h(x,t)+h(x,s)\}
\end{equation}
for all $x\in{\bf R}^N$ and all $0<s<t$.
\end{lemma}
{\bf Proof.}
If $0<s\le t/2$, then 
$t/(t-s)\le 2$ and $s/(t-s)\le 1$, and we obtain 
\begin{eqnarray*}
\Pi_2(t,s:x)
\!\!\! & \le &\!\!\! \frac{t^{\nu/2}}{|t-s|^{\nu/2}} t^{-\frac{\nu}{2}}|G_\alpha(x,t)|
+\frac{s^{\nu/2}}{|t-s|^{\nu/2}} s^{-\frac{\nu}{2}}|G_\alpha(x,s)|\\
\!\!\! & \le &\!\!\! 2^{\frac{\nu}{2}}t^{-\frac{\nu}{2}}|G_\alpha(x,t)|+s^{-\frac{\nu}{2}}|G_\alpha(x,s)|.
\end{eqnarray*}
This together with \eqref{eq:2.6} yields inequality \eqref{eq:7.6} for the case $0<s\le t/2$.
%
%
On the other hand, if $t/2<s<t$, then, 
by the mean value theorem there exists a constant $t_*\in(t/2,t)$ such that
%
%
$$
\Pi_2(t,s:x)
\le|(\partial_t G_\alpha)(x,t_*)|(t-s)^{1-\frac{\nu}{2}}
\le t^{1-\frac{\nu}{2}}|(\partial_tG_\alpha)(x,t_*)|.
$$
This together with \eqref{eq:2.6} implies that
$$
\Pi_2(t,s:x)\le C_1t^{1-\frac{\nu}{2}}
t_*^{-\frac{N}{2}-\frac{2+|\alpha|}{2}}\left[1+(t_*^{-1/2}|x|)^{2+|\alpha|}\right]e^{-\frac{|x|^2}{4t_*}}
\le C_2h(x,t),
$$
for some constants $C_1$ and $C_2$, 
and we obtain inequality \eqref{eq:7.6} for the case $t/2<s<t$. 
Thus Lemma~\ref{Lemma:7.2} follows. 
$\Box$\vspace{5pt}

We are ready to complete the proof of Lemma~\ref{Lemma:2.4}. 
\vspace{5pt}
\newline
{\bf Proof of Lemma~\ref{Lemma:2.4}.}
It suffices to prove \eqref{eq:2.11}. 
We can assume, without loss of generality, that $C_H=1$. 
Let $|\alpha|\le 1$ and 
$$
E(T)=\{(x,y,t,s)\in{\bf R}^{2N}\times(0,T)^2:(x,t)\not=(y,s),\,\, s\le t\}.
$$
By Lemmas~\ref{Lemma:7.1} and \ref{Lemma:7.2} we have 
\begin{eqnarray}
\label{eq:7.7}
\Pi(x,y,t,s) \!\!\! &:= &\!\!\! \frac{|G_\alpha(x,t)-G_\alpha(y,s)|}{|x-y|^\nu+(t-s)^{\nu/2}}
\le\Pi_1(x,y:t)+\Pi_2(t,s:y)\\
 \!\!\! & \preceq &\!\!\! h(x,t)+h(y,t)+h(x,s)+h(y,s)\nonumber
\end{eqnarray}
for all $(x,y,t,s)\in E(T)$. 
On the other hand, by \eqref{eq:2.9} we have  
\begin{eqnarray}
\label{eq:7.8}
 & & \frac{|(\partial_x^\alpha w)(x,t)-(\partial_x^\alpha w)(y,s)|}{|x-y|^\nu+(t-s)^{\nu/2}}\\
 & & \le\int_0^s\int_{{\bf R}^N}
\Pi(x-\xi,y-\xi,t-\tau,s-\tau)H(\xi,\tau)d\xi d\tau\nonumber\\
& & \qquad\qquad
+\int_s^t\int_{{\bf R}^N}
\frac{|G_\alpha(x-\xi,t-\tau)|}{|x-y|^\nu+(t-s)^{\nu/2}}H(\xi,\tau)d\xi d\tau
=:I_1+I_2\nonumber 
\end{eqnarray}
for all $(x,y,t,s)\in E(T)$. 
Then, by \eqref{eq:7.2} and \eqref{eq:7.7} we have 
\begin{eqnarray}
\label{eq:7.9}
 & & 
 I_1\preceq \int_0^s\left(\int_{{\bf R}^N}[h(\xi,t-\tau)+h(\xi,s-\tau)]d\xi\right) d\tau\\
 & & \quad\,
\preceq\int_0^s [(t-\tau)^{-\frac{|\alpha|+\nu}{2}}+(s-\tau)^{-\frac{|\alpha|+\nu}{2}}]d\tau
\preceq s^{1-\frac{|\alpha|+\nu}{2}}
\preceq T^{1-\frac{|\alpha|+\nu}{2}}\nonumber 
\end{eqnarray}
for all $(x,y,t,s)\in E(T)$.  
Furthermore, by \eqref{eq:2.6}
we have
\begin{equation}
\label{eq:7.10}
I_2\preceq\int_s^t
\frac{(t-\tau)^{-|\alpha|/2}}{(t-s)^{\nu/2}}d\tau
\preceq(t-s)^{1-\frac{|\alpha|+\nu}{2}}
\preceq T^{1-\frac{|\alpha|+\nu}{2}}
\end{equation}
for all $(x,y,t,s)\in E(T)$. 
Therefore, by \eqref{eq:7.8}--\eqref{eq:7.10} we have inequality \eqref{eq:2.11}, 
and the proof of Lemma~\ref{Lemma:2.4} is complete.
$\Box$
\bibliographystyle{amsplain}

\end{document}